\documentclass[10pt,a4paper]{amsart}
\usepackage[english, french]{babel}
\usepackage[utf8x]{inputenc}
\usepackage{ucs}
\usepackage{tikz}
\usepackage{amsmath}
\usepackage{amsthm}
\usepackage{amsfonts}
\usepackage{amssymb}
\usepackage{hyperref}
\usepackage{dsfont}
\usepackage{mathrsfs}
\usepackage{eqnarray}
\write18{}
\usepackage{xcolor}
\usepackage{graphicx}
\usepackage{pdfpages} 
\usepackage[nothm]{thmbox}

\hypersetup{
    colorlinks=true,
    linkcolor=red,
    filecolor=magenta,      
    urlcolor=cyan,
}

\numberwithin{equation}{section}

\newcommand{\HH}{\mathbb{H}}

\newcommand{\RR}{\mathbb{R}}

\newcommand{\NN}{\mathbb{N}}
\newcommand{\inte}[1]{\underset{#1}{\int}}

\newcommand{\limisup}[1]{\underset{#1}{\limsup}}
\newcommand{\limiinf}[1]{\underset{#1}{\liminf}}
\newcommand{\infi}[1]{\underset{#1}{\inf}}
\newcommand{\supr}[1]{\underset{#1}{\sup}}
\newcommand{\somme}[1]{\underset{#1}{\sum}}
\newcommand{\limi}[1]{\underset{#1}{\lim}}
\newcommand{\tends}[1]{\underset{#1}{\longrightarrow}}
\newcommand{\equivaut}[1]{\underset{#1}{\sim}}
\newcommand{\fonction}[5]{\begin{array}{r r c l}
					#1 \hspace{2mm} : & #2 & \to & #3 \\
					& #4 & \mapsto & #5 \\
			   \end{array}}
	  
\newcommand{\quotient}[2]{{\raisebox{.2em}{$#1$}\left/\raisebox{-.2em}{$#2$}\right.}}

\DeclareMathOperator{\vol}{vol}
\DeclareMathOperator{\diver}{div}

\theoremstyle{definition}

\newtheorem{definition}[equation]{Definition}
\newtheorem{remark}[equation]{Remark}

\theoremstyle{theorem}

\newtheorem{theoreme}[equation]{Theorem}
\newtheorem{proposition}[equation]{Proposition}
\newtheorem{corollary}[equation]{Corollary}

\newtheorem{lemma}[equation]{Lemma}
\newtheorem*{theoremee}{Theorem}

\newenvironment{poliabstract}[1]
   {\renewcommand{\abstractname}{#1}\begin{abstract}}
   {\end{abstract}}

\title{A stochastic approach to counting problems}

\author{Adrien Boulanger }
\thanks{The author was partially founded by the ERC n°647133 'IChaos'.}
\date{}

\begin{document}

\maketitle

\selectlanguage{english}
\begin{poliabstract}{Abstract} 
We study orbital functions associated to finitely generated geometrically infinite Kleinian groups acting on the hyperbolic space $\mathbb{H}^3$, developing a new method based on the use of the Brownian motion. On the way, we give some estimates of the orbital function associated to nilpotent covers of compact hyperbolic manifolds, partially answering a question asked by M. Pollicott to the author.
\end{poliabstract}

\section{Introduction}
\label{sec introduction}

\textbf{Historical background.} Given a group $\Gamma$ acting properly and discontinuously on a metric space $(X,d)$, we define the \textbf{orbital function} as follows 
	$$ N_{\Gamma}(x,y,\rho) := \sharp \ \{ \ \gamma \in \Gamma \ , \ d(x, \gamma \cdot y) \le \rho \ \} \ , $$ 
where $x,y \in X$ and $\rho > 0$. The growth of orbital functions of groups acting on various types of hyperbolic spaces have been extensively studied since the 50's. Two main approaches were developed in this setting. The oldest one, due to Huber \cite{arthuber} following Delsarte (see \cite{livremartine}), relies on Selberg's pre-trace formula and was designed to study orbital functions in the setting of groups acting cocompactly on the hyperbolic space $\mathbb{H}^2$. This line of work was generalised later on by Selberg (unpublished), Patterson \cite{artpattersondensiteconforme} and Lax-Phillips \cite{articlelax}, among many others. Another way to get estimates of the orbital function growth comes from Margulis' seminal PhD thesis, which exhibits a strong relation between this problem and the mixing of the geodesic flow on the quotient space $\quotient{X}{\Gamma}$. This second approcach has the benefit to give results in the case of variable curvature as well. This idea has been widely extended by numerous authors since, dropping most of Margulis' assumptions required in his work \cite{artyue} \cite{memoireroblin}. For general references, one can recommend Babillot's survey \cite{livremartine} on counting problems, Section 2 of Eskin-McMullen's article \cite{articleeskin} where Margulis' strategy is well explained and the author's PhD dissertion \cite[Chapitre 1]{thesemoi}. The following theorem is, generality-wise, the most advanced of the theory. 

\begin{theoreme}{\cite{memoireroblin}}
\label{theo roblin}
	Let $\Gamma$ be a group acting by isometries, properly and discontinuously on a connected, simply connected complete manifold of negative sectional curvature $X$ such that 
		\begin{itemize}
			\item the length spectrum of the quotient space is non arithmetic;
			\item the unitary tangent bundle of the quotient space $\quotient{X}{\Gamma}$ admits a finite Bowen-Margulis-Sullivan measure.
		\end{itemize}
Then, for any points $x,y \in X$ there exists a constant $C$ such that
	$$ N_{\Gamma}(x,y,\rho) \equivaut{\rho \to \infty} C \ e^{ \delta_{\Gamma} \rho} \ ,$$
where $\delta_{\Gamma}$ is the critical exponent of the group $\Gamma$. 
\end{theoreme}

Recall that the critical exponent is, in full generality, defined as 
$$ \delta_{\Gamma} := \limisup{\rho \to \infty} \ \frac{\ln(N_{\Gamma}(x,y,\rho))}{\rho}  \ ,$$ 
which was proven to be a limit in \cite{artroblinfoncorb}. The above theorem is weaker than the one stated in \cite{memoireroblin} but the author did not want to burden this article with definitions regarding the CAT(-1) setting for which an analogous version holds. Let us note that the two assumptions of Theorem \ref{theo roblin} are automatically satisfied when $\Gamma$ acts convex-compactly on a hyperbolic space of any dimension, see \cite[Proposition 3]{articleguivarch} for the assumption made on the non-arithmeticity of the length spectrum. \\

Roblin also showed, under the same assumptions of those of Theorem \ref{theo roblin}, that if the quotient space $\quotient{X}{\Gamma}$ does not admit any finite Bowen-Margulis-Sullivan measure, one has 
	\begin{equation*}
		 N_{\Gamma}(x,y,\rho) = o \left( e^{ \delta_{\Gamma} \rho} \right) \ ,
	\end{equation*}	 
for any pair of points $x,y \in X$ when $\rho \to \infty$. Nowadays, the goal is to explicit the asymptotic in this case. There is, up to the author's knowledge, only two classes of examples where the finitness assumption had been successfully dropped. The first one is due to Pollicott and Sharp \cite{artsharppolliabeliancover}. They found an asymptotic for orbital functions of groups associated to abelian covers of compact hyperbolic manifolds. The second one is due to Vidotto \cite{thesevidotto} and deals with an example in variable curvature. In both cases, the authors manage to get an asymptotic of the orbital function. Their methods rely on a finite measure property hidden somewhere: in \cite{artsharppolliabeliancover} the author's method uses a strong mixing property of the geodesic flow of the compact underlying manifold whereas in \cite{thesevidotto} the geodesic flow still admits a coding by a sub-shift of finite type. \\

\textbf{The strategy and the setting.} In this article we introduce a method to study the growth of orbital functions of groups acting on the real hyperbolic space $\mathbb{H}^3$ of dimension three, based on the use of the Brownian motion in the spirit of Sullivan's work \cite{articlesullivanpositivity}. Our motivation for introducing the Brownian motion in this setting is to compensate the lack of known randomness of the geodesic flow - preventing to adapt Margulis' method - when there is no finite Bowen-Margulis-Sullivan measure. Replacing the geodesic flow by the Brownian motion gives a larger range of applications, the later being random on its own. The key facts allowing us to relate its dynamical behaviour to the orbital function is the \textbf{strong homogeneity} of the hyperbolic space and the so called \textbf{drift property} of the Brownian motion. We refer to the author's PhD dissertation \cite{thesemoi} for more details. \\

More precisely, this article aims at investigating the link bewteen Brownian motion and orbital functions through examples. We call \textbf{Kleinian group} any torsion free subgroup $\Gamma$ of the orientation preserving isometries of the hyperbolic 3-space $\mathbb{H}^3$ acting properly and discontinuously or, equivalently, the fundamental group of a hyperbolic manifold of dimension 3. Regarding the counting problem for finitely generated Fuchsian groups, analoguous to the last but acting on $\mathbb{H}^2$, there is not so much to be said since all the quotients $\quotient{\mathbb{H}^2}{\Gamma}$ turned out to carry a finite Bowen-Margulis-Sullivan measure \cite{dalhorgeo}, so that Theorem \ref{theo roblin} readily applies in this case. Therefore, one would like to picture what is the situation in dimension 3. We will call a \textbf{finitely generated} Kleinian group $\Gamma$ \textbf{degenerate} when the quotient manifold $ M_{\Gamma} := \quotient{\mathbb{H}^3}{\Gamma}$ does not carry any finite Bowen-Margulis-Sullivan measure. Note that we require $\Gamma$ to be finitely generated in order to qualify it degenerate. By extension, we also call degenerate the manifold $M_{\Gamma}$. These manifolds are central objects in some of the most famous theorems around the theory of finitely generated Kleinian groups, such as Thurston's hyperbolisation theorem (see \cite{livremcmumu} or \cite{livreotal2001hyperbolization}), Ahlfors' conjecture, corollary \cite{articlecanary} of Marden's tameness conjecture (see \cite{surveycanarymarden} and the reference therein) or the ending lamination conjecture \cite{articleendinglamination}, all theorems now. \\

Surprisingly, very little is known about the behaviour of their orbital functions except that they have critical exponent $2$ \cite{articlesullivandemicylindre}, which is to say that its exponential growing rate is of the order the volume growth of $\HH^3$. The only special case already settled being as a corollary of \cite{artsharppolliabeliancover}, giving an asymptotic to the orbital function for surface groups associated to $\mathbb{Z}$-covers of compact 3-manifolds, included in M. Pollicott and R. Sharp's theorem as Abelian covers of compact manifolds. \\

\textbf{Description of the results.} Our two main theorems can be stated as follows. We refer to Subsection \ref{subsec their geometry} for the definition and the description of ends of hyperbolic manifolds involved in the statements.

\begin{theoreme}
	\label{theo kleinian counting 1}
	Let $\Gamma$ be a degenerate Kleinian group such that the manifold $M_{\Gamma}$ has positive injectivity radius and has both degenerate and geometrically finite ends. Then, for any $x,y \in \mathbb{H}^3$ there are constant $C > 0$ and $ \rho_0 > 0$ such that for any $\rho > \rho_0$ we have
	$$ N_{\Gamma}(x,y,\rho) \le C \ \frac{e^{2 \rho}}{\rho} \ .$$
\end{theoreme}

Under some extra assumption, one can improve the above theorem into

\begin{theoreme}
	\label{theo kleinian counting 2}
Let $\Gamma$ be a degenerate Kleinian group such that the manifold $M_{\Gamma}$ has positive injectivity radius, and the \textbf{averaged orbital function is roughly decreasing} (see Subsection \ref{subsec sharper estimate}), then 
\begin{itemize}
	\item either all the ends of $M_{\Gamma}$ are degenerate, then there are constants $C_-, C_+$ such that for any $x,y \in \mathbb{H}^3$ there is $\rho_0 > 0$ such that for $\rho > \rho_0$ we have
		$$ \frac{C_- \ e^{2\rho} }{\sqrt{\rho}} \le N_{\Gamma}(x,y,\rho) \le \frac{C_+ \ e^{2\rho} }{\sqrt{\rho}} ; $$
	\item or, for every $x,y \in M_{\Gamma}$ there are constants $C_-, C_+, \rho_0 > 0$ such that for $\rho > \rho_0$ we have
			$$ \frac{C_- \ e^{2\rho} }{\rho^{\frac{3}{2}}} \le N_{\Gamma}(x,y,\rho) \le \frac{C_+ \ e^{2\rho} }{\rho^{\frac{3}{2}}} \ . $$
\end{itemize}	 

\end{theoreme}

Note that the constants appearing in the first above item depend only on $\Gamma$. Constants involved in the second item depend on the points under consideration, which is not a flaw of the method but something expected: we shall come back on why through the proof of Theorem \ref{theo chaleur Kleinian} (see, in particular, Figure \ref{figuremixedtype}). \\

The rough decreasing assumption is natural since there is no known orbital function not satisfying it. However, the author does not know how to prove that an orbital function is roughly decreasing without using a counting theorem, which is disappointing. \\

To prove both Theorem \ref{theo kleinian counting 1} and Theorem \ref{theo kleinian counting 2} one needs first to link the heat kernel to the orbital function, which is the purpose of the next two theorems that we will prove in Section \ref{sec orbital function heat kernel}.

\begin{theoreme}
		\label{theo upper bound}
	There is a constant $C_+ > 0$ such that for any Kleinian group $\Gamma$ and all $x,y \in \mathbb{H}^3$ one has for $\rho > 1$
				\begin{equation}
						\label{eqtheo upper bound}
						\frac{N_{\Gamma}(x,y,\rho)}{e^{2 \rho}} \le C_+ \sqrt{\rho} \ p_{\Gamma} \left( x, y, \frac{\rho}{2} \right)  \ ,
				\end{equation} 
		where $p_{\Gamma}(x,y,t)$ is the heat kernel of the manifold $M_{\Gamma} = \quotient{\mathbb{H}^3}{\Gamma}$.
\end{theoreme}

\begin{remark} 
This upper bound is relevant only when the critical exponent of the group is $2$. Otherwise, it would not even have the correct exponential growth since from \cite[Theorem 2.17]{articlesullivanpositivity} and \cite[Theorem 10.24]{bookgrigorheatkernel} we have
$$ \frac{\ln p_{\Gamma}(x,y,\rho)}{\rho} \tends{\rho \to \infty} \delta_{\Gamma} (\delta_{\Gamma} -2) \ . $$
 One might, perhaps, get an more general upper bound in weighting the manifold $M_{\Gamma}$ with respect to its first eigenfunction, known to be positive \cite[Theorem 2.1]{articlesullivanpositivity}.
\end{remark}

The above upper bound is not sharp, and finer estimates - sharp upper and lower bounds up to multiplicative constants - can be proven under a polynomial control of the heat kernel in large time and some extra assumption on the orbital function. We refer to Subsection \ref{subsec sharper estimate} for the definition involved in second item of the theorem below.

\begin{theoreme}
\label{theo refined estimates}
Let $\Gamma$ be a Kleinian group. If we have both
\begin{itemize}
	\item there is $\alpha \ge 0$ such that there is two points $x,y \in M_{\Gamma}$, two constants $c_-, c_+$ and $t_0 > 0$ such that for $t > t_0$ one has 					
		\begin{equation}
			\label{eq bound heat kernel}
			 c_- \ t^{-\alpha} \le p_{\Gamma}(x,y,t) \le c_+ \ t^{-\alpha} \ ;
		 \end{equation}
	\item the averaged orbital function is roughly decreasing,
\end{itemize}
then there are positive constants $C_-,C_+, \rho_0$ such that for any $\rho > \rho_0$ we have
	$$ \frac{C_- \ e^{2\rho}}{\rho^{\alpha}} \le N_{\Gamma}(x,y,\rho) \le  \frac{C_+ \ e^{2\rho}}{\rho^{\alpha}} \ . $$
\end{theoreme}

Ultimately, the proof of Theorem \ref{theo upper bound} and Theorem \ref{theo refined estimates} relies on the drift property of the Brownian motion and on the fact that $\HH^3$ is locally symmetric. Note that we did not require $\Gamma$ to be finitely generated in the assumptions of Theorem \ref{theo upper bound} and Theorem \ref{theo refined estimates}. \\

Both Theorem \ref{theo upper bound} and Theorem \ref{theo refined estimates} apply immediately in the context of a manifold admitting the volume doubling property and a Poincaré inequality (see Subsection \ref{subsecpoincaredoubling}), which is the case for  nilpotent covers of compact manifold. 

\begin{corollary}
	\label{cor counting nilpotent}
		Let $\Gamma_0$ be a Kleinian group acting co-compactly and $\Gamma \triangleleft \Gamma_0$ a normal subgroup of $\Gamma_0$ such that $N := \quotient{\Gamma_0}{\Gamma}$ is nilpotent. Then, there is $\alpha \in \NN$ (only depending on $N$) and $C_+ >0$ such that for any $\rho > 1$ 
		$$ N_{\Gamma}(x,y,\rho) \le C_+ \ \frac{\ e^{2\rho}}{ \rho^{\frac{\alpha -1}{2}}} \ . $$
Moreover, if the averaged orbital function is roughly decreasing, then there exists two constant $C_-, C_+$ such that for every $x,y$ there is ${\rho}_0$ such that for $\rho > \rho_0$ we have 
		$$  \frac{C_- \ e^{2 \rho}}{\rho^{\alpha/2}} \  \le N_{\Gamma}(x,y,\rho) \le \frac{C_+ \ e^{2 \rho}}{\rho^{\alpha/2}} \ . $$
\end{corollary}

In order to deduce our main theorems \ref{theo kleinian counting 1} and \ref{theo kleinian counting 2} from Theorems \ref{theo upper bound} and \ref{theo refined estimates}, we investigate the large time behaviour of the heat kernel on the manifolds under consideration, the degenerate hyperbolic manifolds. The proof of the following theorem will occupy all the section \ref{sec proofs}.

\begin{theoreme}
\label{theo chaleur Kleinian}
Let $\Gamma$ a degenerate Kleinian group such that $M_{\Gamma}$ has positive injectivity radius, then 

\begin{itemize}
	\item either all the ends of $M_{\Gamma}$ are degenerate, then there is two constant $C_-, C_+ > 0$ such that for every $x,y \in \mathbb{H}^3$ there is $t_0 > 0$ such that for $t > t_0$ 
		$$ \frac{C_- }{\sqrt{t}} \le p_{\Gamma}(x,y,t) \le \frac{C_+}{\sqrt{t}} \ ; $$
	\item or there are degenerate and geometrically finite ends, then for every $x,y \in M_{\Gamma}$ there exists a constant $C_+ > 0$ such that one has for $t > 0$ 
			$$ p_{\Gamma}(x,y,t) \le \frac{C_+}{t^{\frac{3}{2}}} \ . $$
		Moreover, for any $x \in M_{\Gamma}$ there is a constant $ C_- > 0$ such that for $t > 0$ 
		$$ \frac{C_-}{t^{\frac{3}{2}}} \le p_{\Gamma}(x,x,t) $$
\end{itemize}	 
\end{theoreme}

In order to show the above theorem, we rely on the knowledge of the geometry of the ends and on the fact that we can compare $M_{\Gamma}$ to simple graphs for which we know these bounds (or their analytical counterpart) to hold. \\

If $M_{\Gamma}$ has only degenerated ends it must be roughly isometric to a graph consisting of a base point on which we attached finitely many copies of $\NN$. We show that this graph satisfies both the Poincaré inequality and the doubling property, see Subsection \ref{subsecpoincaredoubling}. Both these properties are invariant under rough isometries and therefore hold for the manifold $M_{\Gamma}$ itself. We conclude by using a well known theorem giving a very precise two sided estimate for the heat kernel in this setting. \\

The case where there is both degenerated and geometrically finite ends is more subtle. In fact, most of the theorems giving an upper bound for the heat kernel actually give an uniform upper bound, which is not expected in our case, see Figure \ref{figuremixedtype}. We circumvent this difficulty in weighting the hyperbolic manifold with respect to a harmonic function constructed by Sullivan and Thurston. The associate weighted heat kernel and the usual one are precisely related by a time independent factor, which allows one to compare their large time behaviour. We conclude by showing  that an uniform upper bound actually holds for the weighted manifold in comparing it to a graph. We recover the lower bound by using the upper one. \\

\textbf{Acknowledgement.} I am especially grateful to the anonymous referee who simplified and corrected many proofs in the first version of this article. I want to thank G. Carron and L. Saloff-Coste for helping me navigating the theory of heat kernel on manifolds. I would like to thank my PhD advisors, G. Courtois and F. Dal'Bo, as well as S. Gouëzel for numerous comments and suggestions about this article. I finally want to thank F. Ledrappier for stimulating conversations around the exponent $3/2$. 

\section{Weighted manifolds and their heat kernels.}

This section is about to introduce an important object for what concerns this article; heat kernels of weighted Riemannian manifolds. 

\label{sec weighted heat kernel}

\subsection{Basics on weighted heat kernels.} We denote by $(M,g)$ a complete Riemannian manifold and by $\mu_g$ its associated Riemannian measure. Given a positive function $\sigma$ on $M$, the \textbf{weight}, we call the metric measured space $(M,g, \sigma \mu_g)$ a \textbf{weighted Riemannian manifold}. We denote the weighted Riemannian measure $\sigma \mu_g$ by $\mu_{\sigma}$ for short. The following definition aims at enlarging the Laplace operator to the setting of weighted Riemannian manifolds. We denote by $L^2(M,\mu_{\sigma})$ the space of square integrable functions of $M$ with respect to the measure $\mu_{\sigma}$. 

\begin{definition}
	The \textbf{weighted Laplace operator}, denoted by $\Delta_{\sigma}$, is defined as follows
		$$ \Delta_{\sigma} f := - \sigma^{-1} \diver ( \sigma \nabla f)$$
	where $f$ is a smooth function on $M$.
\end{definition}

We refer to the non-weighted (equivalently, weighted with the constant function 1) Riemannian Laplace operator as the \textbf{usual Laplace operator}. \\

The following theorem leads to the definition of the heat kernel in the setting of a weighted manifold. 

\begin{theoreme}{\cite[Corollary 4.11, page 117]{bookgrigorheatkernel}}
\label{theocauchyprobleml2}
	Let $(M,g,\mu_{\sigma})$ be a complete weighted Riemannian manifold. For any initial condition $u_0 \in L^2(M, \mu_{\sigma})$ the Cauchy problem
	\begin{equation}
		\label{eq Cauchy problem}
		 \left\{ 
			\begin{array}{l}
		 		\Delta_{\sigma} u(x,t) + \partial_t u(x,t) = 0 \\
		 		u(\cdot,t) \tends{t \to 0} u_0 \hspace{0.5 cm} \text{ in \ }  \ L^2 \ . 
			\end{array} 
		 \right. 
	\end{equation}
	
has a unique solution $u : M \times \RR_+ \to \RR$ given by 
 $$u(t,x) := \inte{M} \ p_{\sigma}(x,y,t) \ f(y) \ d \mu_{\sigma}(y) \ , $$
where $p_{\sigma}(x,y,t)$ is, by definition, the \textbf{heat kernel} of the weighted manifold $(M,g, \mu_{\sigma})$. 
\end{theoreme}

As a comprehensive reference on heat kernels one can recommend \cite[Chapters 6 and 7]{bookgrigorheatkernel}. 
\\

The usual heat kernel, meaning the one coming from the above theorem but for the usual Laplace operator, would be simply denoted by $p(x,y,t)$. \\

Provided that the weights of two manifolds differ from one another by a positive harmonic function, the following well known proposition relates the two induced heat kernels by a time independent factor. In particular, $x$ and $y$ being fixed, the large time behaviour of these heat kernels are alike. 

\begin{proposition}
	\label{proposition heat kernels related}
	Let $(M,g)$ be a complete Riemannian manifold and $h$ a positive harmonic function. The weighted heat kernel $p_{h^2}(x,y,t)$ associated to $(M,g,\mu_{h^2})$ satisfies 
		\begin{equation}
			\label{eqprop weighted heat kernel}
			 p_{h^2}(x,y,t) = \frac{p(x,y,t)}{h(x)h(y)} \ .
		\end{equation}
\end{proposition}	

The proof relies on uniqueness in Theorem \ref{theocauchyprobleml2}: one must verify that the above kernel satisfies the Cauchy problem \ref{eq Cauchy problem}, which is straightforward and left to the reader.

\subsection{Rough isometries.} \label{subsec rough isometries} The following equivalence relation is of crucial importance for us since two metric measured spaces equivalent under this relation would turn out to share a lot of geometrico-analytic properties simultaneously. This would allow us to prove these geometrico-analytic properties for simpler model equivalent to the weighted manifold we want to work with. In our setting, these models would turn out to be weighted graphs. If $(X,d)$ is a metric space, we denote by $B(x, r)$ the ball centred at $x \in X$ of radius $r > 0$.

\begin{definition}[Kanaï \cite{artkanai}]
	Let $(M_1,d_1,\mu_1)$ and $(M_2,d_2,\mu_2)$ two  metric measured spaces. We will say of an application $f : M_1 \to M_2$ that it is a \textbf{rough isometry} if the three following conditions hold:
		\begin{enumerate}
			\item there is $\epsilon > 0$ such that the $\epsilon$-neighbourhood of $f(M_1)$ is equal to $M_2$;
			\item there is $a,b > 0$ such that for every $x,y \in M_1$ we have
				$$ a^{-1} \ d_1(x,y) -b \le d_2(f(x),f(y)) \le a \ d_1(x,y) + b \ . $$
			\item there is $C > 0$ such that for every $x \in M_1$ we have $$ C^{-1} \mu_1 \Big( B(x, 1) \Big) \le \mu_2 \Big( B(f(x), 1) \Big) \le C \ \mu_1 \Big( B(x, 1) \Big) . $$ 
		\end{enumerate}
\end{definition}

This definition is very close to the more familiar one of quasi-isometry which only requires (1) and (2). The notion of rough isometry may be considered as an extension of the one of quasi-isometry to the setting of metric measured spaces. \\

This notion was also developed in order to understand the large time behaviour of the heat kernel thanks to the behaviour of the transition kernel of a random walk on a graphs coarsely looking alike the manifold. There is a natural construction of such a graph, called a discretisation, see \cite[Section 1]{articlegillesdiscretisation}. Recall how to construct it. Given a weighted Riemannian manifold $(M,g,\mu_{\sigma})$ we construct a weighted graph $(G,m)$ as follows: given an $\epsilon >0$ we define the vertices of $G$ as any $\epsilon$-separated maximal family of points of $M$. An edge relates two vertices if and only if the two points under consideration are at distance at most $2 \epsilon$. The weight $m$ is defined a
	$$ m(x) := \mu_{\sigma}(B(x_M, \epsilon)) \ ,$$
where $x_M$ is the point on $M$ corresponding to the vertex $x \in G$. A discretisation, as a graph, comes with a structure of metric measured space in taking the graph distance and the weighted counting measure $ m := \somme{x \in G} m(x) \delta_x$.

\begin{lemma}[Lemma 2.5 \cite{artkanai}]
	A discretisation $G$ of a Riemannian manifold $(M,g)$ with Ricci curvature bounded below is roughly isometric to $(M,g)$.
\end{lemma}

\subsection{Poincaré inequality and doubling property.}
\label{subsecpoincaredoubling} This section gathers the results we use in order to prove the part of Theorem \ref{theo chaleur Kleinian} which addresses the case where all ends are degenerate. \\

In \cite{artcoulhonsaloffisometrie} the authors proved that many large scale properties of a metric measured space are invariant under rough isometries, like the doubling volume property and the Poincaré inequality defined below. 

\begin{definition}
	A graph $G$ satisfies the \textbf{doubling property} if there is $C > 0$ such that for all $x \in G$ and all $r > 0$ we have
	\begin{equation}
			\label{eq doubling} 
				\sharp \ B(x,2 r) \le C \ \sharp \ B(x, r)  \ .
		\end{equation}
\end{definition}

In order to state what is a Poincare inequality, we define the discrete analogous of the derivative  for a graph
	$$ \delta f (x) :=  \sqrt{\somme{y \sim x} \ |f(y) - f(x)|^2} $$
where $x \sim y$ stands for the relation '$x$ and $y$ are related by an edge of the graph'. 

\begin{definition}
A graph $G$ satisfies a \textbf{Poincaré inequality} if there is $P > 0$ such that for any $x \in G$, $r > 0$ and for any function $f$ on $G$ we have
	\begin{equation}
				\label{eqpoincare}
					 \somme{y \in B(x,r)} (f(y) - f_r)^2 \le P \  r^2  \somme{y \in B(x, 2 r)} (\delta f)(y)^2  \ ,
	\end{equation}
			where $f_r$ is the mean of $f$ on the ball $B(x,r)$.
\end{definition}

The following stated theorem is a combination of two theorems. The first one \cite{artcoulhonsaloffisometrie} asserts that if a Riemannian manifold with Ricci curvature bounded from below is roughly isometric to a graph which satisfies the Poincaré inequality and the doubling property then it satisfies the parabolic Harnack inequality. This property is known to be equivalent to a very precise control of the heat kernel by \cite{articlegrigonilpotent}, see also \cite[Lecture II]{surveysaloffcoste}. In particular, it implies the two sided estimate appearing in the following theorem, which is enough for our purpose.

\begin{theoreme}{ \cite{articlegrigonilpotent}, \cite[Theorem 4.2]{articlesaloffpoincaresobolev}, \cite[Theorem 8.2]{artcoulhonsaloffisometrie}}
\label{theo grigo equivalence}
	Let $(M,g)$ be a Riemannian manifold with Ricci curvature bounded below which is roughly isometric to a graph which satisfies both the doubling property \eqref{eq doubling} and the Poincare inequality \eqref{eqpoincare}. Then, there are constants $C_-, C_+ > 0$ such that for any $x,y \in M$ there is $t_0 > 0$ such that for any $t > t_0$ we have
	$$  \frac{C_- }{\mu_{g} \Big( B(x, \sqrt{t}) \Big)}  \le p(x,y,t) \le  \frac{C_+ }{\mu_{g} \Big( B(x, \sqrt{t}) \Big)} \ . $$
\end{theoreme} 

In the setting of a cover of a compact manifold, the covering manifold is roughly isometric to any Cayley graph $\mathcal{C}(N)$ of the deck group $N$. In our case the deck group $N$ is nilpotent. The volume growth of nilpotent groups are very well understood \cite[Section 1]{artbaas}: there are constant $c_-, c_+$ and an integer $\alpha$ depending only on the nilpotent covering deck group $N$ such that 
	$$c_- \rho^{\alpha/2} \le \vol \Big( B_{\mathcal{C}(N)}(\sqrt{\rho}) \Big) \le c_+ \rho^{\alpha/2} \ ,$$
where $\vol \Big( B_{\mathcal{C}(N)}(\rho) \Big)$ is the volume of any ball of radius $\rho$ in the graph $\mathcal{C}(N)$. In particular one can immediately deduce from this that the Cayley graph satisfies both the doubling property and the Poincare inequality using \cite[Theorem 2.2]{artkleiner}. Therefore, one recover Corollary \ref{cor counting nilpotent} by using the theorems \ref{theo upper bound} and \ref{theo refined estimates} and the fact that polynomial volume growth is invariant under rough isometries \cite[Proposition 2.2]{artcoulhonsaloffisometrie}. \\

We will also use Theorem \ref{theo grigo equivalence} in Section \ref{sec proofs} to handle the part of Theorem \ref{theo chaleur Kleinian} which addresses the case where all the ends are degenerate.

\subsection{Heat kernel of weighted manifolds.} This section gathers the results we use in order to prove the part of Theorem \ref{theo chaleur Kleinian} which addresses the case where there is both degenerated ends and geometrically finite ones. We will ultimately rely on the comparison theorem \ref{theocomparison} stated below. This theorem would allow us to transfer the desired upper bound for the case of interest from the geometric properties of a graph to which it is roughly isometric. \\

Before stating it, we need to introduce two 'semi-local' properties analogous to \ref{eq doubling} and \ref{eqpoincare} that the manifold must carry in order to this theorem.

\begin{definition}
	A weighted Riemannian manifold $(M,g, \mu_{\sigma})$ is said to satisfy 
	\begin{itemize} 
		\item the $(DV)_0$ \textbf{property} if for any $r_0 > 0$ there is $C = C(r_0) > 0$ such that for all $x \in G$ and all $r < r_0$ we have 
	\begin{equation}
			\label{eqlocaldoubling} 
				 \mu_{\sigma}( B(x,2 r)) \le C \ \mu_{\sigma}( B(x, r))  \ 
		\end{equation}
	\item the $(P_2)_0$ \textbf{property} if for all $r_0 >0$ there is a constant $ P =P(r_0) > 0$ such that for any $x \in M$, any $ 0 <r <r_0$ and for any compactly supported and smooth function $f$ 
	\begin{equation}
				\label{eqlocalpoincare}
					 \inte{B(x,r)} (f(y) - f_r(x))^2 d \mu_{\sigma} \le P \  r^2  \inte{ B(x, 2 r)} | \nabla f|^2  d \mu_{\sigma} \ ,
	\end{equation}
			where $f_r(x)$ is the mean of $f$ on the ball $B(x,r)$ taken with respect to $\mu_{\sigma}$. 
	\end{itemize}
\end{definition}

Theorem \ref{theocomparison} comes also with an assumption on the graph to which we want to compare our manifold to. This assumption is the 'non local' part of the theorem. Given a measured space $(M,\mu)$ and a function $f : M \to \RR$ we denote for any $p \ge 0$ 
	$$ || f ||_{L^p(\mu)} := \left( \inte{M} |f|^p \ d \mu \right)^{\frac{1}{p}} \ .$$

\begin{definition}
	A weighted graph $(G,m)$ is said to carry \textbf{the Sobolev inequality} $S_{q,p}$ if there is a constant $C = C(q,p) > 0$ such that for any compactly supported function $f$ on $G$ we have 
		$$ ||f||_{L^p(m)} \le C || \delta f||_{L^q(m)} \ .$$
\end{definition}	

Let us now state the comparison theorem. Note that the following theorem is not stated in the generality of weighted Riemannian manifolds in \cite{artcoulhonsaloffisometrie}. However, it also holds in this case since all the results of \cite{artcoulhonsaloffisometrie} can be generalised to the setting of roughly isometric operators see \cite[Section 9]{artcoulhonsaloffisometrie} (which is even more general that the setting of weighted Riemannian manifolds). 

\begin{theoreme}{\cite[Corollaire 5.6, Proposition 6.5]{artcoulhonsaloffisometrie}}
	\label{theocomparison}
Let $(M,g,\mu_{\sigma})$ be a weighted manifold satisfying both $(DV)_0$ and $(P_2)_0$ which is roughly isometric to a weighted graph $G$. then for any $\nu > 2$ the two following properties are equivalent.
\begin{itemize}	
	\item there are $t_0, C$ such that for any $t > t_0$ we have
	$$ \underset{x, y \in M}{\sup} p_{\sigma}(x,y,t) \le C \ t^{- \nu/2} \ ; $$ 
   \item the weighted graph $(G,m)$ carries a Sobolev inequality $S_{2,p}$ with $p = \frac{2\nu}{\nu-2}$.
 \end{itemize}
\end{theoreme}
		
The above theorem would be the key to recover the upper bound of Theorem \ref{theo chaleur Kleinian} in the case where there is both geometrically finite ends and degenerate one. For what concerns the lower bound, we will use

\begin{theoreme}{\cite{articlegrigocoulhon} \cite[Theorem 16.6 page 427]{bookgrigorheatkernel}}
\label{theo lower bound}
	Let $(M,g,\mu_\sigma)$ be a complete weighted manifold and $x \in M$ such that 
\begin{itemize}
	\item the weighted manifold satisfies the doubling property at $x$, \textit{i.e.} there is a constant $C > 0$ such that for every $r > 0$ 
		$$ \mu_{\sigma}(B(x,2r)) \le C \mu_{\sigma}(B(x,r)) \ ;$$
	\item there is a constant $C_+ > 0$ such that for every $t > 0$
		$$p_{\sigma}(x,x,t) \le \frac{C_+}{\mu_{\sigma} \Big(B(x,\sqrt{t}) \Big) } \ , $$
\end{itemize}

then there exists a constant $C_-$ such that for any $t > 0$ one has 
		$$p_{\sigma}(x,x,t) \ge \frac{C_-}{\mu_{\sigma} \Big( B(x,\sqrt{t}) \Big)} \ . $$
\end{theoreme}

\section{The orbital function and the usual heat kernel} 

\label{sec orbital function heat kernel}

We start off reviewing some properties of the heat kernel on $\mathbb{H}^3$ and its quotient manifolds; hyperbolic manifolds of dimension 3. We will then prove successively the upper bound of Theorem \ref{theo upper bound} and theorem \ref{theo refined estimates}, controlling the orbital function provided some knowledge of the large time behaviour of the heat kernel. On the way, we will define the averaged orbital function and the rough decreasing property required to apply Theorem \ref{theo refined estimates}. These are ones of the main two ingredients to prove Theorems \ref{theo kleinian counting 1} and \ref{theo kleinian counting 2}. 

\subsection{Heat kernels of hyperbolic manifolds.}
The hyperbolic 3-space $\HH^3$ enjoys an explicit formula for its heat kernel denoted by $p_3(x,y,t)$.

\begin{theoreme}[see \cite{artgrigoheatkernelhyper}]
	\label{theo heat kernel 3-hyperbolic}
Let $\Gamma$ be a Kleinian group.	For any two points $x,y \in \HH^{3}$ at distance $\rho$ and for all $t > 0$ one has
			\begin{equation}
				\label{eqtheo heat kernel 3-hyperbolic}
					p_{3}(x,y,t) = \frac{1}{ (4 \pi t)^{\frac{3}{2}}} \frac{\rho}{\sinh(\rho)} e^{ -  t - \frac{\rho^2}{4t}} \ .
			\end{equation}
\end{theoreme}

For $\rho \ge 0$ and $t > 0$, we denote by $p_3(\rho,t)$ the function giving the value of the heat kernel $p_3(x,y,t)$ for any two points $x,y \in \mathbb{H}^3$ at distance $\rho$ from one another. 

\begin{remark}
Note that heat kernels of hyperbolic spaces of different dimensions are related by the following formula (see \cite{artgrigoheatkernelhyper}):
\begin{equation}
	\label{eqformule reccurence noyau}
		 p_{n+2}(\rho,t) = - \frac{ e^{ -nt }} {2 \pi \sinh(\rho) } \partial_{\rho} p_{n}(\rho,t) \ . 
\end{equation}
There is also an explicit formula for the heat kernel of the hyperbolic space of dimension 2, so that, using the above formula one can derive a formula in any dimension. The results of this article can therefore be extended to any dimension. However, our main application being toward finitely generated Kleinian groups, we will stick to the dimension three in order not to burden this article with quite a lot of computations.
\end{remark}

\begin{lemma}
	\label{lemma heat kernel quotient}
Let $\Gamma$ be a Kleinian group. For any two points $x,y$ on $M_{\Gamma}$ the heat kernel of the manifold $M_{\Gamma}$ is given by the formula:
		\begin{equation}
			\label{eqlemma heat kernel quotient}
			 p_{\Gamma}(x,y,t) := \somme{\gamma \in \Gamma} \ p_{3}(\tilde{x}, \gamma \cdot \tilde{y},t) \ .  
		\end{equation}
where $\tilde{x}$ and $\tilde{y}$ are any two lifts on $\mathbb{H}^3$ of the points $x$ and $y$.
\end{lemma} 

\textbf{Proof.} The proof is straightforward. One first needs to show that the above series makes sense, which follows from explicit formula given by \ref{theo heat kernel 3-hyperbolic} which shows that the heat kernel decreases spatially super-exponentially whereas the orbital function $N_{\Gamma}(x,y,\rho)$ as at most exponential growth in the variable $\rho$. Then one needs to show that the formula does not depend on the choices of the lifts $\tilde{x}$ and $\tilde{y}$ which follows, for the variable $\tilde{y}$, from the fact that we averaged with respect to $\Gamma$. For the variable $\tilde{x}$ it follows from the fact that $\Gamma$ acts by isometries. One can verify that the above series satisfies the Cauchy problem \eqref{eq Cauchy problem} using the fact that $p_3(x,y,t)$ satisfies the equivalent Cauchy problem on $\HH^3$. We conclude by uniqueness.   \hfill $\blacksquare$ 

\subsection{The orbital function upper bound}

This subsection is entirely dedicated to prove the upper bound of Theorem \ref{theo upper bound} announced in the introduction, that we recall here for the reader's convenience.

\begin{theoremee} There is a universal constant $C_+ > 0$ such that for any Kleinian group $\Gamma$ and for any $x,y \in \mathbb{H}^3$ one has for $t>1$
				\begin{equation*}
						\frac{N_{\Gamma}(x,y,t)}{e^{2 t}} \le C_+ \sqrt{t} \ p_{M_\Gamma} \left( x, y, \frac{t}{2} \right) \ .
				\end{equation*}
\end{theoremee}

\textbf{Proof.} We start in proving the following

\begin{lemma}
	\label{lemmelinkintegral}
 Under the same assumptions as in the above theorem, we have for any $x,y \in M_{\Gamma}$ and any $t > 0$ 
\begin{equation}
	\label{equationlinkintegral}
	 p_{\Gamma}(x, y, t)  = - \inte{\mathbb{R}_+} N_{\Gamma}(\rho,x,y) \ \partial_{\rho} p_3 
			(\rho,t) d \rho \ . 
\end{equation}	
\end{lemma}		

\textbf{Proof.} Since the group $\Gamma$ acts discretely, the function 
$$ g := t \mapsto N_{\Gamma}(x,y,t)$$ 
has bounded variation. Therefore, one can apply the Stieltjes integral by part formula to any smooth function $f$ integrated with respect to $dg$ to get 

	\begin{equation*}
		\begin{split}
			 \somme{ \gamma \in \mathcal{A}_{\Gamma}(x,y,a,b) } f \left( d(x, \gamma \cdot y) \right) = & \\ 
		 N_{\Gamma}(x,y,b) f(b) - & N_{\Gamma}(x,y,a) f(a)  - \inte{[a,b]} N_{\Gamma}(x,y,\rho) \ f'(\rho) d \rho  \ .
		\end{split}
	\end{equation*}

We specialise the above identity with the data 
$$ \fonction{f}{\mathbb{R}_+}{\mathbb{R}}{\rho}{p_3(\rho,t) \ ,} $$
and $b \to + \infty$, which gives 
\begin{equation*}
			 \somme{ \gamma \in \Gamma } \ p_3 \left( x, \gamma \cdot y, t \right) = 
		 \limi{b \to + \infty} \left( N_{\Gamma}(x,y,b) p_3(b,t)  - \inte{[0,b]} N_{\Gamma}(x,y,\rho) \ \partial_{\rho} p_3(\rho,t) d \rho \right) \ .
\end{equation*}

The heat kernel $p_3(x,y,\rho)$ being decreasing spatially super-exponentially, one also has 
	 $$\limi{b \to + \infty} N_{\Gamma}(x,y,b) p_3(b,t) = 0 \ ,$$
	
and therefore
\begin{equation*}
			p_\Gamma \left( x, y, t \right) = - \limi{b \to + \infty} \inte{[0,b]} N_{\Gamma}(x,y,\rho) \ \partial_{\rho} p_3(\rho,t) d \rho \ .
\end{equation*}

From Equation \ref{eqformule reccurence noyau} we have
		$$ - N_{\Gamma}(x,y,\rho) \partial_{\rho} p_3(\rho,t) \ge 0 \ ,$$ 
	
so that, letting $b \to \infty$, we get 
	\begin{equation*}
			 	p_{\Gamma}(x, y, t)  = - \inte{\RR_+} N_{\Gamma}(x,y,\rho) \ \partial_{\rho} p_3 (\rho,t) d \rho \ ,
	\end{equation*}
which is the desired result. \hfill $\blacksquare$  \\

In order to deduce from Lemma \ref{lemmelinkintegral} Theorem \ref{theo upper bound} we use the explicit formula for the heat kernel on $\HH^3$ as announced. The key feature, from our perspective, given by Formula \eqref{eqtheo heat kernel 3-hyperbolic} is the drift of the Brownian motion. Brownian paths tend to concentrate in large time on an annulus of small radius $2t - \sqrt{t}$ and of large radius $2t + \sqrt{t}$. This $\sqrt{t}$-lack-of-concentration on the perfect sphere of radius $2t$ is responsible for the non optimality of the upper bound of Theorem \ref{theo upper bound}. We shall come back later on this and see how to strengthen the result adding an extra assumption on the orbital function. \\

Because the function $-\partial_\rho p_3(\rho,t)$ is positive, one has the following lower bound of the left member of \eqref{equationlinkintegral}:
\begin{equation}
	\label{eqprooupboundroughmin}	
	 p_{\Gamma}(x, y,t) \ge   \inte{[2 t, 2 t +1]} N_{\Gamma}(x,y,\rho) \left[ - \partial_{\rho} p_3(\rho,t) \right] d \rho \ .
\end{equation}

The orbital function $\rho \mapsto N_{\Gamma}(x,y,\rho)$ being increasing and non negative, one also has
		$$  p_{\Gamma}(x, y,t) \ge  N_{\Gamma}(x,y, 2 t)  \inte{[2 t , 2 t+1]} -\partial_{\rho}p_3(\rho,t) d \rho \ . $$ 
	
Therefore, the existence of a constant $C_1 >0$ such that for every $\rho \in [2t,2t +1]$ 
		\begin{equation} 
			\label{eqprooftheo upper bound}
				\left[ -\partial_{\rho} p_3 \right]  ( \rho, t ) \ge \frac{C_1}{\sqrt{t}} e^{ - 4 t } \ ,
		\end{equation}

would give 
$$  p_{\Gamma}(x, y,t) \ge  N_{\Gamma}(x,y, 2 t)  \frac{C_1}{\sqrt{t}} e^{ - 4 t } \ , $$ 

and thus the upper bound \eqref{eqtheo upper bound} of Theorem \ref{theo upper bound} in substituting $t$ with $ 2 t$. \\

To show that the lower bound \eqref{eqprooftheo upper bound} holds, we recall Formula \ref{eqtheo heat kernel 3-hyperbolic}:
		$$p_{3}(\rho,t) = \frac{1}{ (4 \pi t)^{\frac{3}{2}}} \frac{\rho}{\sinh(\rho)} e^{ -  t - \frac{\rho^2}{4t}} \ , $$
		
hence, 
\begin{align*}
	\partial_{\rho}p_{3}(\rho,t) = \frac{1}{ (4 \pi t)^{\frac{3}{2}}} \left( \frac{1}{\sinh(\rho)} - \frac{\rho \cosh(\rho)}{\sinh(\rho)^2} - \frac{\rho^2}{2 t \sinh(\rho)} \right) e^{ -  t - \frac{\rho^2}{4t}} 
\end{align*}	

and then,
\begin{equation}
	\label{equation derive heat}
	\partial_{\rho}p_{3}(\rho,t) = \frac{1}{ (4 \pi t)^{\frac{3}{2}}} \left( 1 - \frac{\rho \cosh(\rho)}{\sinh(\rho)} - \frac{\rho^2}{2 t} \right) \frac{e^{ -  t - \frac{\rho^2}{4t}}}{\sinh (\rho)} \ .
\end{equation}

Note that there are constants $C_2, C_3 > 0$ such that for every $t >1$ and $\rho \in [2t, 2t+1]$ we have:
 $$ \frac{\rho \cosh(\rho)}{\sinh(\rho)} + \frac{\rho^2}{2 t}  -1 \le C_2 \left( 2t + \frac{4t^2}{2 t} \right) \le C_3 t \ , $$
 
from which we get a constant $C_4 > 0$ such that for every $t >1$ and $\rho \in [2t, 2t+1]$ we have:
\begin{align*}
	- \partial_{\rho}p_{3}(\rho,t) & \ge \frac{ 2 C_3  t}{ (4 \pi t)^{\frac{3}{2}}} e^{-2t} e^{ -  t - \frac{4t^2}{4t}} \\
			 & \ge \frac{C_4}{ t^{\frac{1}{2}}} e^{ -4 t } \ ,
\end{align*}	 
which concludes the proof.  \hfill $\blacksquare$ \\ 

As already emphasised, this theorem may be re-enforced under some extra assumptions.

\subsection{Sharper estimates toward finer results}

\label{subsec sharper estimate}

Under a 'rough decreasing' assumption defined below, one can sharpen the proof of Theorem \ref{theo upper bound} to get a better control of the orbital function growth. 
\begin{definition}
A function $f : \RR_+ \to \mathbb{R}_+$ is said to be \textbf{roughly decreasing} if there are positive constants $C > 0$ and $ \rho_0 > 0$ such that for any $\rho_0 \le \rho_1 \le \rho_2$ one has 
		$$ f(\rho_2) \le C \ f(\rho_1) \ .$$
\end{definition}

This rough decreasing property will be required for the \textbf{averaged orbital function} defined as
	$$ \tilde{N}(x,y,\rho) := \rho \mapsto \frac{N(x,y,\rho)}{e^{2\rho}} \ . $$

Recall that $e^{2 \rho}$ is roughly the volume of a hyperbolic ball of radius $\rho$. Before entering the proof of Theorem \ref{theo refined estimates}, we recall here its statement for the reader's convenience.

\begin{theoremee}
Let $\Gamma$ be a Kleinian group such that there is $\alpha \ge 0$ such that there is two points $x,y \in M_{\Gamma}$ such that
\begin{itemize}
	\item  there are constants $c_-, c_+, t_0$ such that for $t \ge t_0$
		\begin{equation}
			\label{eq bound heat kernel}
			 c_- \ t^{-\alpha} \le p_{\Gamma}(x,y,t) \le c_+ \ t^{-\alpha} \ ;
		 \end{equation}
	\item the averaged orbital function $ \rho \mapsto \tilde{N}(x,y,\rho)$ is roughly decreasing,
\end{itemize}
then there is two positive constants $C_-$ and $C_+$ such that 
	$$ \frac{C_- \ e^{2\rho}}{\rho^{\alpha}} \le N_{\Gamma}(x,y,\rho) \le  \frac{C_+ \ e^{2\rho}}{\rho^{\alpha}} \ . $$
\end{theoremee}

We shall see that the first condition made on the heat kernel is fulfilled in natural examples, as nilpotent covers of a given compact hyperbolic manifold or for thick degenerate Kleinian groups. \\

\textbf{Proof.} We prove first the upper bound appearing in the conclusion of Theorem \ref{theo refined estimates} and then we will use it to prove the lower one. For what follows the points $x,y$ are fixed, so that all the functions appearing along the proof have to be seen as functions defined on $\mathbb{R}_+$. Given two functions $f,g : \mathbb{R}_+ \to \mathbb{R}_+$, we denote $f(t) \prec g(t)$ the transitive relation : there are constant $\rho_0, C>0$ such that $ f(t) \le C g(t)$ for $t \ge t_0$, omitting the value of the constants. \\

The proof of the upper bound is rather similar to the proof of Theorem \ref{theo upper bound}. Recall first Equation \ref{equationlinkintegral}:
	\begin{equation*}
			 	p_{\Gamma}(x, y, t)  = - \inte{\RR_+} N_{\Gamma}(x,y,\rho) \ \partial_{\rho} p_3 (\rho,t) d \rho \ .
	\end{equation*}
	
Using the upper bound assumed on the heat kernel, one has 
	\begin{equation*}
		 \inte{\RR_+} N_{\Gamma}(x,y,\rho) \ [-\partial_{\rho} p_3 (\rho,t)] d \rho \prec \frac{1}{t^{\alpha}} \ .
	\end{equation*}

Instead of the rough lower bound giving by Equation \eqref{eqprooupboundroughmin} we prefer 
	\begin{equation*}
		 \inte{2t - \sqrt{t}}^{2t+\sqrt{t}} N_{\Gamma}(x,y,\rho) \ [-\partial_{\rho} p_3] (\rho,t) d \rho \prec \frac{1}{t^{\alpha}} \ ,
	\end{equation*}

that can be rewritten as
	\begin{equation*}
		 \inte{2t - \sqrt{t}}^{2t+\sqrt{t}}  \tilde{N}_{\Gamma}(x,y,\rho)  \ e^{2 \rho} [-\partial_{\rho} p_3] (\rho,t) d \rho \prec \frac{1}{t^{\alpha}} \ .
	\end{equation*}

Using the rough decreasing assumption made on the averaged orbital function we get 
 \begin{equation}
 	\label{equ prop finer 1}
		\tilde{N}_{\Gamma}(x,y,2t + \sqrt{t}) \inte{2t - \sqrt{t}}^{2t+\sqrt{t}} - e^{2 \rho} \partial_{\rho} p_3 (\rho,t) d \rho \prec \frac{1}{t^{\alpha}} \ . 
 \end{equation}
	
Recalling Equation \eqref{equation derive heat} one has
\begin{align*}
	e^{2 \rho} \ \partial_{\rho}p_{3}(\rho,t) & = \frac{1}{ (4 \pi t)^{\frac{3}{2}}} \left( 1 - \frac{\rho \cosh(\rho)}{\sinh(\rho)} - \frac{\rho^2}{2 t} \right) \frac{e^{ -  t - \frac{\rho^2}{4t} + 2 \rho}}{\sinh (\rho)} \ ,
\end{align*}	

from which we deduce the existence of two universal constants $C_1, C_2 > 0$ such that for every $t >1$ and $\rho \in (t , 3t)$ 
 \begin{equation}
 	\label{equ prop finer 2}
	  \frac{ C_1 }{ t^{\frac{1}{2}}} e^{ -  t - \frac{\rho^2}{4t} + \rho}  \le - e^{2 \rho} \partial_{\rho}p_{3}(\rho,t)  \le \frac{C_2}{ t^{\frac{1}{2}}} e^{ -  t - \frac{\rho^2}{4t} + \rho} \ . 
\end{equation}

In particular,
  $$ \inte{2t - \sqrt{t}}^{2t+\sqrt{t}} - e^{2 \rho} \partial_{\rho} p_3 (\rho,t) d \rho \succ 1  \ .$$
 
So that Equation \eqref{equ prop finer 1} gives
$$ \tilde{N}_{\Gamma}(x,y,2t + \sqrt{t}) \prec \frac{1}{t^{\alpha}} \ , $$

which implies 
$$ \tilde{N}_{\Gamma}(x,y,t) \prec \frac{1}{t^{\alpha}} \ , $$

since 
\begin{align}
	\label{eqlimsup}
	 \limisup{ t \to \infty} \ \tilde{N}_{\Gamma}(x,y,t) \ t^{\alpha} & = \limisup{ t \to \infty} \ \tilde{N}_{\Gamma}(x,y,2t+ \sqrt{t}) \ (2t + \sqrt{t})^{\alpha}  \\
	 & \le 	 \limisup{ t \to \infty} \ N_{\Gamma}(x,y,2t+ \sqrt{t}) (2t)^{\alpha} \ . 
\end{align}

The strategy to get the lower bound is rather similar to what we have done for the upper one, but more subtle. Let us mention that we will need the already proven upper bound in order to get the lower one. \\

We cannot bound from below as roughly as we did for the upper bound, we prefer to start off the mid-rough approximation 
	\begin{equation*}
			 	p_{\Gamma}(x, y, t)  \equivaut{ t \to \infty} - \inte{t}^{3t} N_{\Gamma}(x,y,\rho) \ \partial_{\rho} p_3 (\rho,t) d \rho \ ,
	\end{equation*}
	
valid since the we assumed that the heat kernel has polynomial decay and because the two remaining integrals are decreasing faster than any rational function (from the very rough upper bound $N_{\Gamma}(x,y,\rho) \prec e^{2\rho}$). From the polynomial lower bound assumed on the heat kernel, we also have
$$    \frac{1}{t^{\alpha}}  \prec  \inte{t}^{3t} N_{\Gamma}(x,y,\rho) \ [-\partial_{\rho} p_3] (\rho,t) d \rho \ . $$

We now use the right side of Inequality \eqref{equ prop finer 2} to get 
$$    \frac{1}{t^{\alpha}} \prec \inte{t}^{3t} \frac{N_{\Gamma}(x,y,\rho)}{e^{2 \rho}} \ \frac{ e^{ -  t - \frac{\rho^2}{4t} + \rho} }{ t^{\frac{1}{2}}}  \ d \rho \ , $$

which can be rewritten  
\begin{equation}
	\label{eqpropfiner3}
    \frac{1}{t^{\alpha}} \prec \inte{t}^{3t} \tilde{N}_{\Gamma}(x,y,\rho) \ \frac{ e^{ -  t - \frac{\rho^2}{4t} + \rho} }{ t^{\frac{1}{2}}}  \ d \rho \ . 
\end{equation}

Let chop off the above integral the following way, $k>0$ will be chosen later on.
	
	\begin{align*}
		 \inte{t}^{3t} \frac{N_{\Gamma}(x,y,\rho)}{e^{2\rho}}  \frac{ e^{ -  t - \frac{\rho^2}{4t} + \rho}  }{ t^{\frac{1}{2}}} \ d \rho & =  \left. \inte{t}^{2t - k \sqrt{t}} \frac{N_{\Gamma}(x,y,\rho)}{e^{2\rho}}  \frac{ e^{ -  t - \frac{\rho^2}{4t} + \rho}  }{ t^{\frac{1}{2}}} \ d \rho \ \right\} := I_1(t) \\
		 & + \left. \inte{2t - k\sqrt{t}}^{2t + k \sqrt{t}} \frac{N_{\Gamma}(x,y,\rho)}{e^{2\rho}}  \frac{ e^{ -  t - \frac{\rho^2}{4t} + \rho} }{ t^{\frac{1}{2}}} \ d \rho \ \ \right\} := I_2(t)  \\
		& +  \left. \inte{2t + k \sqrt{t}}^{3t} \frac{N_{\Gamma}(x,y,\rho)}{e^{2\rho}}  \frac{ e^{ -  t - \frac{\rho^2}{4t} + \rho}}{ t^{\frac{1}{2}}} \ d \rho \ \ \right\} := I_3(t) \ .
		 \end{align*}
		 
The goal is now to show the
	\begin{lemma} 
 		\label{lemmeproofsharper}
		$$I_i(t) \prec \frac{e^{-k^2/4}}{t^{\alpha}}$$
for $i=1,3$.
	\end{lemma} 
	
Before proving the above lemma, let us show how to conclude using almost the same argumentation than the one used to get the upper bound. Combined with Equation \eqref{eqpropfiner3} the above lemma gives 
	$$ \frac{1}{t^{\alpha}} \prec \frac{e^{-k^2}}{t^{\alpha}} + I_2(t) + \frac{e^{-k^2}}{t^{\alpha}} \ , $$
taking $k$ large enough we get $ t^{-\alpha} \prec I_2(t)$, namely:
	$$ \frac{1}{ t^{\alpha}} \prec \inte{2t - k\sqrt{t}}^{2t + k \sqrt{t}} \frac{N_{\Gamma}(x,y,\rho)}{e^{2\rho}}  \frac{ e^{ -  t - \frac{\rho^2}{4t} + \rho}  }{ t^{\frac{1}{2}}} \ d \rho \ . $$

Using the rough decreasing assumption made one has 
\begin{align*}
	\frac{1}{t^{\alpha}} \prec \tilde{N}_{\Gamma}(x,y,2t -k \sqrt{t})  \inte{2t - k\sqrt{t}}^{2t + k \sqrt{t}} \frac{  e^{ -  t - \frac{\rho^2}{4t} + \rho} }{ t^{\frac{1}{2}}} \ d \rho \ .
 \end{align*}  

The right above integral being bounded, we also have 
$$ \frac{1}{t^{\alpha}}  \prec \tilde{N}_{\Gamma}(x,y,2t -k \sqrt{t}) \ .  $$

We conclude, the same way as done for the upper bound with Equations \eqref{eqlimsup} but using the liminf instead of the limisup. \hfill $\blacksquare$ \\

It remains to show that Lemma \ref{lemmeproofsharper} holds. \\

\textbf{Proof. of of Lemma \ref{lemmeproofsharper}.} The proofs for $I_1(t)$ and $I_3(t)$ are very similar: we only perform the computation concerning $I_1(t)$. We start off using the already proven upper bound $N_{\Gamma}(x,y,\rho) \prec \frac{ e^{2 \rho}}{\rho^{\alpha}}$:
	\begin{align*}
		I_1(t) & :=  \inte{t}^{2t - k \sqrt{t}} \frac{N_{\Gamma}(x,y,\rho)}{e^{2\rho}}  \frac{ e^{ -  t - \frac{\rho^2}{4t} + \rho} }{ t^{\frac{1}{2}}} \ d \rho \\ 
			& \prec \inte{t}^{2t - k \sqrt{t}} \frac{1}{\rho^{\alpha}} \frac{e^{ - (\frac{2t - \rho}{2\sqrt{t}})^2}}{\sqrt{t}} \ d \rho \\
			& \prec \frac{1}{t^{\alpha}} \ \inte{t}^{2t - k \sqrt{t}} \frac{e^{ - (\frac{2t - \rho}{2\sqrt{t}})^2}}{\sqrt{t}} \ d \rho \ .
	\end{align*}
We substitute $\rho$ by $u := \frac{2t - \rho}{2\sqrt{t}}$ in the above integral to get
\begin{align*}
\inte{t}^{2t - k \sqrt{t}} \frac{e^{ - (\frac{2t - \rho}{2\sqrt{t}})^2}}{\sqrt{t}} d \rho & = \inte{k/2}^{\sqrt{t}} \frac{e^{ -u^2 }}{\sqrt{t}} \ 2 \sqrt{t} \ d u  \\
    & \le 2 \inte{k/2}^{\infty} e^{ -u^2 } \ d u  \\
	& \le 2 \ e^{-k^2/4} \ ,
\end{align*}	

which gives  $$ I_1(t) \prec \frac{e^{-k^2/4}}{t^{\alpha}} \ ,$$

the expected result. \hfill $\blacksquare$ $\blacksquare$

\section{The thick tamed degenerate hyperbolic manifolds.}
\label{sec thick tamed}

This section is devoted to the description of the topology and the geometry of degenerate hyperbolic 3-manifolds with positive injectivity radius, called \textbf{thick}. Note in particular that it prevents parabolic elements. It basically gathers results from the literature. 

\subsection{Their geometry.} 
\label{subsec their geometry} The celebrated tameness theorem, a result proved by Agol and Calegari-Wise independently asserts that the quotient $M_{\Gamma}$ of $\mathbb{H}^3$ by a finitely generated Kleinian group $\Gamma$ is homeomorphic to the interior of a compact manifold with boundary. One can recommend \cite{surveycanarymarden} for a detailed presentation of this theorem and its history. On the topological side, this theorem asserts that if $\Gamma$ is a finitely generated Kleinian group then $M_{\Gamma}$ can be described as a compact submanifold $K$ with boundary on which one has attached finitely many ends $E_1, ..., E_p$ all of them homeomorphic to $S \times \RR_+$ where $S$ is a surface. \\

 On the geometric side, these ends can be put in two families depending on whether or not they belong to the convex core of $M_{\Gamma}$, convex core which is constructed as follows. Given $\Gamma$ a Kleinian group acting on $\HH^3$ we define the limit set $\Lambda_{\Gamma}$ of $\Gamma$ as the the set in $\mathbb{S}^2 = \partial_{\infty} \HH^3$ on which accumulates any orbit $\Gamma \cdot x_0$ where $x_0$ is any point of $\HH^3$ (it is easy to see that the limit set does not depend on $x_0$). Let $C(\Gamma)$ be the smallest convex set of $\HH^3$ which contains all the infinite geodesics of $\HH^3$ having endpoints in $\Lambda_{\Gamma}$: the convex core $CC(M_{\Gamma})$ is then defined as $\quotient{C(\Gamma)}{\Gamma}$ (which is well defined since $\Gamma$ preserves $\Lambda_{\Gamma}$ by construction and, therefore, $C(\Gamma)$). One can recommend Thurston's note [Chapter 8], \cite[Part II.1]{livrecanary_marden_epstein_2006} and \cite[Appendix A]{artmoroianuschlenker} for a detailed presentation of convex sets in $\HH^3$, convex cores and their properties. An end contained in the convex core is called \textbf{degenerate} and an end which is not contained in it is called \textbf{geometrically finite}, see Figure \ref{figdefends}. We introduce the following terminology.
\begin{itemize}
 \item if all the ends of $M_{\Gamma}$ are degenerate, we will say that $\Gamma$ is \textbf{fully degenerate};
 \item if $M_{\Gamma}$ carries both degenerate ends and geometrically finite ones, we will say that $\Gamma$ is of \textbf{mixed type}. 
\end{itemize}

Note that the convex core of a manifold associated to a fully degenerate Kleinian group is the whole manifold. Therefore, we could also have defined alternatively the dichotomy fully degenerate/mixed type as whether or not the convex core is the entire manifold, which is equivalent to whether or not the limit set of the Kleinian group is the entire sphere. This last alternative definition, involving the limit set, is perhaps the most standard. There is a lot of references on ends of hyperbolic manifold with finitely generated groups. One can recommend the introduction of \cite{articlecanary}, the book \cite{livreotal2001hyperbolization} or Mahan Mj' survey. \\

\begin{figure}[h!]
\begin{center}
	\def\svgwidth{0.6 \columnwidth}
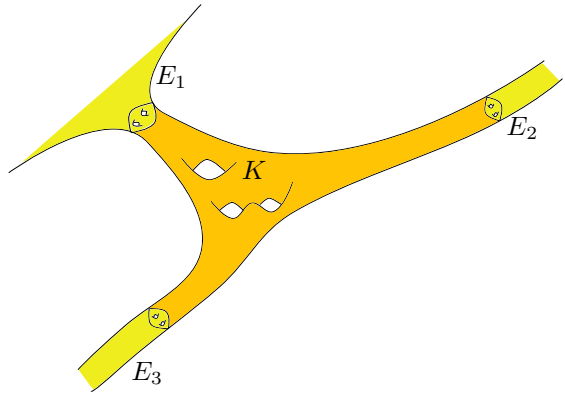
	\end{center}
\caption{This manifold carries 3 ends, $E_1$, $E_2$ and $E_3$. The first one is geometrically finite and the two last ones are degenerate. The orange part, corresponding to a compact submanifold on which we have glued the ends, is represented in orange. The convex core $CC(M_{\Gamma})$ consists of the union of $K$ together with the degenerate ends $E_2$ and $E_3$ or, equivalently, to the complement of the end $E_1$.}
	\label{figdefends}
\end{figure}

Let us start in describing more precisely the geometrically finite ends which are well more understood than the degenerate ones. They basically features the same geometry than the 2-dimensional funnels. \\

For our purpose, we will only need to know that they have positive bottom of the spectrum. Recall that for a given Riemannian manifold with boundary $(M,g)$, the bottom of the spectrum is defined as 
$$	 \lambda_0(M) := \infi{ f \in \mathcal{C}_0^{\infty}(M)} \frac{|| \nabla f||^2_2}{|| f||^2_2} \ , $$ 

where $\mathcal{C}_0^{\infty}(M)$ stands for the set of smooth and compactly supported function vanishing on the boundary. The terminology comes from the fact that the above defined number is the actual bottom of the spectrum of the Friedrichs extension of the Laplace operator with Dirichlet boundary conditions, see for example \cite[Chapter 7]{livrerosenberg}. \\

If $A \subset \HH^3$ and $C \ge 0$ we denote by $A_{C}$ the $C$ neighbourhood of the set $A$;
	$$ A_C := \{  \ x \in \HH^3 \ , \ d(x,A) \le C \ \} \ . $$

\begin{proposition}
	\label{propspectralgap}
	The hyperbolic manifold with boundary $M_{\Gamma} \setminus CC(M_{\Gamma})_1$ has positive bottom of the spectrum.
\end{proposition}

In particular any subset of $M_{\Gamma} \setminus CC(M_{\Gamma})_1$ has positive bottom of the spectrum and, therefore, all the geometrically finite ends. \\

\textbf{Proof.} The proof relies on the well known fact that if there is a smooth positive function $f_{\Gamma} : \HH^3 \setminus CC(M_\Gamma)_1 \to \RR^*_+$ which satisfies 
\begin{equation}
	\label{equation bottomspectrum}
	 \Delta f_{\Gamma} \ge c f_{\Gamma} 
\end{equation}
for some $c > 0$ then $\lambda_0(M_{\Gamma} \setminus CC(M_{\Gamma})_1) \ge c$, see for example \cite[Proposition 3.2]{arttapieroblin}. \\

In our setting, we will construct such a function by smoothing the distance function to the convex core $CC(\Gamma)$. If $A \subset \HH^3$ we denote by $d_A :  \HH^3 \to \RR_+$ the function distance to $A$.

\begin{lemma}
	\label{lemma trouspectral}
	Let $K$ be a convex set of $\HH^3$. There is $\alpha$ small enough such that the function
		\begin{equation}
			\label{equation trouspectral}
		 e^{-\alpha d_K} : \HH^3 \setminus K_{1/2} \to \RR_+^* 
		\end{equation}
satisfies \eqref{equation bottomspectrum} for some $c > 0$ in the distributional sense.
\end{lemma}

Before proving the above lemma, let us see how to use it in order to deduce that Proposition \ref{propspectralgap} holds. Since $C(\Gamma)$ is a convex set by construction, one would be able to conclude provided that the function $e^{-\alpha d_{C(\Gamma)}}$ is both smooth and invariant under the action of the group $\Gamma$. \\

The difficulty is that the function $d_{C(\Gamma)}$ may not be smooth since the boundary of a convex set is not smooth in general. Let us then first see how to smooth the function $s := e^{- \alpha d_{C(\Gamma)}}$. \\

Let $\kappa : (0,1/2) \to \RR_+$ be any compactly supported smooth function. We define  
	$$ \tilde{s}(x) := \inte{\HH^3} \kappa(d(x,y)) \ s(y) \ d \mu_g(y) \ ,$$
which is a smooth function from $\HH^3 \to \RR_+$. Since $\HH^3$ is symmetric one has 
	$$ \Delta_x \ \kappa(d(x,y)) = \Delta_y \ \kappa(d(x,y)) \ . $$
In particular for any $x \in \HH^3 \setminus C(\Gamma)_1$ we have
\begin{align*}
	 \Delta(\tilde{s})(x) & =  \inte{\HH^3} \Delta_y \kappa(d(x,y)) \ s(y) \ d \mu_g(y) \\
	 & \ge c \inte{\HH^3}  \kappa(d(x,y)) \ s(y) \ d \mu_g(y) =  c \ \tilde{s}(x) \ ,
\end{align*}
since $\kappa$ is supported on $(0, 1/2)$ and because \eqref{equation trouspectral} holds for $s$ on $ \HH^3 \setminus C(\Gamma)_{1/2}$ (in the distributional sense). Therefore 
	$$ \fonction{f}{\HH^3 \setminus CC(M_{\Gamma})_1}{\RR_+^*}{x}{\tilde{s}(x)} $$
is smooth and satisfies \eqref{equation bottomspectrum}. Let us now see that $f$ is $\Gamma$-invariant. \\

Note that $s$ was $\Gamma$-invariant since $C(\Gamma)$ is $\Gamma$-invariant by construction (the limit set $\Lambda_{\Gamma}$ is itself invariant). Note also that, because $\kappa(d(\gamma \cdot x, y)) = \kappa( d(x, \gamma^{-1} \cdot y))$, the map $\tilde{s}$ inherits the $\Gamma$-invariance from the $\Gamma$-invariance of $s$. \\

The function $f$ (which is the restriction of $\tilde{s}$ to $\HH^3 \setminus CC(M_{\Gamma})_1$) is then $\Gamma$-invariant and therefore descends to a smooth positive function $ f_{\Gamma} : M_{\Gamma} \setminus CC(M_{\Gamma})_1 \to \RR_+$ which also satisfies \eqref{equation bottomspectrum},  concluding. It remains then to prove Lemma \ref{lemma trouspectral}. \\

\textbf{Proof of Lemma \ref{lemma trouspectral}.} We compute 
	$$ \Delta ( e^{-\alpha d_K} ) = ( - \alpha \Delta d_K - \alpha^2 |\nabla d_K|^2)  e^{-\alpha d_K} \ . $$

As the distance to a given set, the function $d_K$ satisfies $|\nabla d_K|^2 = 1$. Therefore 
	$$ \Delta ( e^{-\alpha d_K} ) = ( - \alpha \Delta d_K - \alpha^2)  e^{-\alpha d_K} \ .$$
One will be able to conclude by taking $\alpha$ small enough provided that there is $c > 0$ such that for any $x \in \HH^3 \setminus  K_{1/2}$ one has (in the distributional sense)
$$- \Delta d_K \ge c \ ,$$
which is what we will focus on. \\ 

We adapt and detail an argument used in \cite[Lemma A.12]{artmoroianuschlenker} to the setting of interest. The remark made in this article is that the Hessian of the distance function to a convex set can be compared to the Hessian of the distance function to some totally geodesic hyperplane. \\

Let $x$ be any point in $\HH^3 \setminus K$ and $p$ be its nearest point projection on the convex set $K$ (which is well known to be well defined, see for example \cite[Lemma A11]{artmoroianuschlenker}). Let $\nu$ be the unique geodesic from $x$ to $p$ and $\mathcal{P}$ the hyperplane orthogonal to the geodesic $\nu$ which passes by the point $p$ as in the following figure.

\begin{figure}[h!]
\begin{center}
	\def\svgwidth{0.6 \columnwidth}
\begingroup%
  \makeatletter%
  \providecommand\color[2][]{%
    \errmessage{(Inkscape) Color is used for the text in Inkscape, but the package 'color.sty' is not loaded}%
    \renewcommand\color[2][]{}%
  }%
  \providecommand\transparent[1]{%
    \errmessage{(Inkscape) Transparency is used (non-zero) for the text in Inkscape, but the package 'transparent.sty' is not loaded}%
    \renewcommand\transparent[1]{}%
  }%
  \providecommand\rotatebox[2]{#2}%
  \newcommand*\fsize{\dimexpr\f@size pt\relax}%
  \newcommand*\lineheight[1]{\fontsize{\fsize}{#1\fsize}\selectfont}%
  \ifx\svgwidth\undefined%
    \setlength{\unitlength}{634.48627142bp}%
    \ifx\svgscale\undefined%
      \relax%
    \else%
      \setlength{\unitlength}{\unitlength * \real{\svgscale}}%
    \fi%
  \else%
    \setlength{\unitlength}{\svgwidth}%
  \fi%
  \global\let\svgwidth\undefined%
  \global\let\svgscale\undefined%
  \makeatother%
  \begin{picture}(1,0.66899428)%
    \lineheight{1}%
    \setlength\tabcolsep{0pt}%
    \put(0,0){\includegraphics[width=\unitlength,page=1]{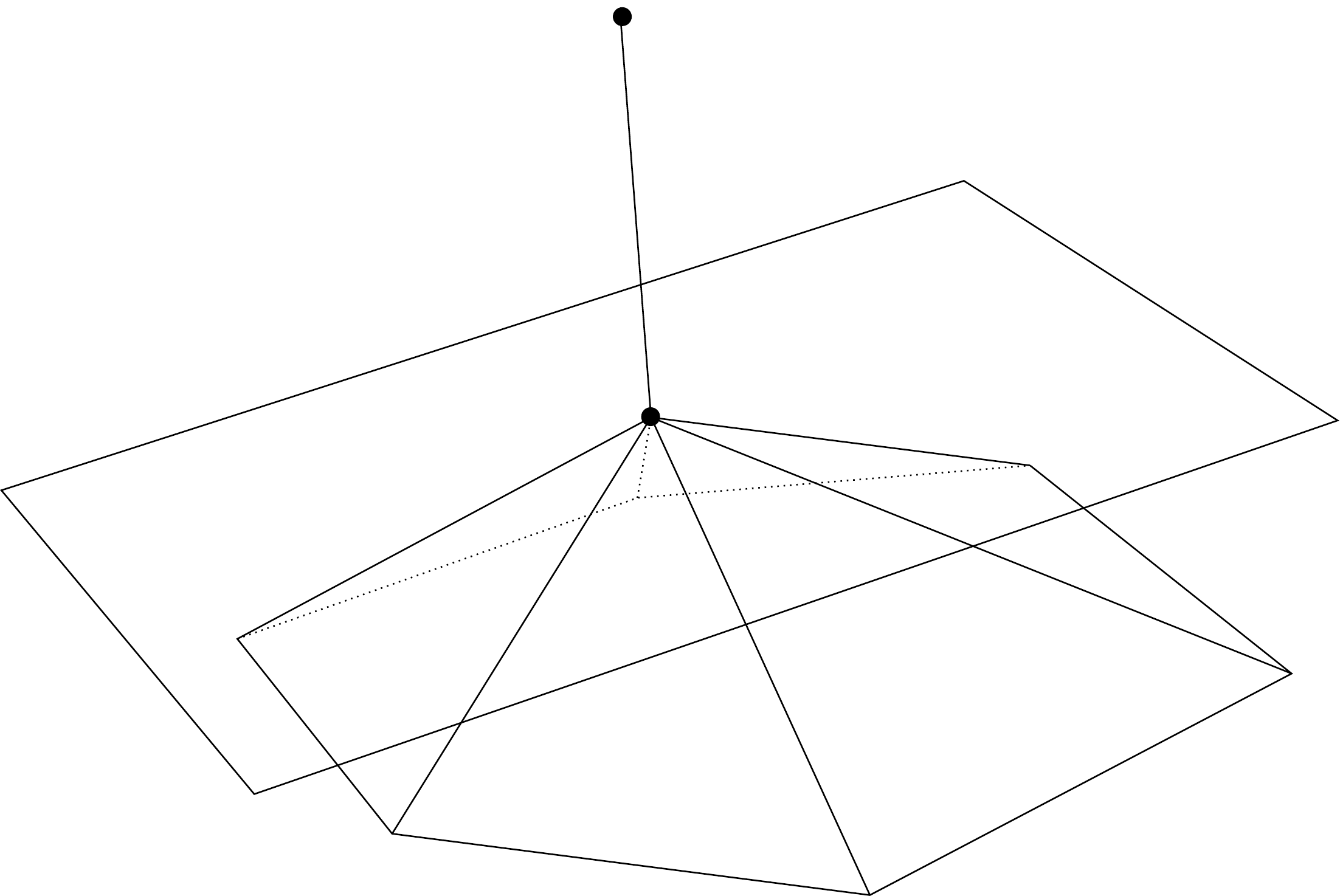}}%
    \put(0.4871524,0.63260916){\color[rgb]{0,0,0}\makebox(0,0)[lt]{\lineheight{1.25}\smash{\begin{tabular}[t]{l}$x$\end{tabular}}}}%
    \put(0.50828891,0.37192636){\color[rgb]{0,0,0}\makebox(0,0)[lt]{\lineheight{1.25}\smash{\begin{tabular}[t]{l}$p$\end{tabular}}}}%
    \put(0.08415086,0.2775169){\color[rgb]{0,0,0}\makebox(0,0)[lt]{\lineheight{1.25}\smash{\begin{tabular}[t]{l}$\mathcal{P}$\end{tabular}}}}%
    \put(0.68583494,0.13237995){\color[rgb]{0,0,0}\makebox(0,0)[lt]{\lineheight{1.25}\smash{\begin{tabular}[t]{l}$K$\end{tabular}}}}%
    \put(0.48760497,0.53397245){\color[rgb]{0,0,0}\makebox(0,0)[lt]{\lineheight{1.25}\smash{\begin{tabular}[t]{l}$\nu$\end{tabular}}}}%
  \end{picture}%
\endgroup%

	\end{center}
	\label{figdefends}
\end{figure}

Note that, by convexity of $K$, one has for all $y$ in the half space containing $x$ whose boundary is $\mathcal{P}$ 
\begin{equation}
	\label{eqtrouspectralconvex}
	 d_{\mathcal{P}}(y) \le d_K(y) 
\end{equation}
with equality in $x$.  We denote by $c(t)$ any geodesic such that $c(0) = x$.  For any set $A \subset \HH^3$ and any curve $c$ we also denote 
	$$ \Phi_c^A(t) := d_A(c(t)) \ .$$
Inequality \eqref{eqtrouspectralconvex} together with the equality case implies in particular that for any geodesic $c$ such that $c(0) = x$ we have (in the distributional sense) 
	$$ (\Phi_c^{K})''(0) \ge (\Phi_c^{\mathcal{P}})''(0) \ .$$
Indeed, the function $ f := \Phi_c^{K} - \Phi_c^{\mathcal{P}} : \RR \to \RR$ is non negative and satisfies $f(0) = 0$ which implies $f''(0) \ge 0$ (in the distributional sense). \\

We get in particular 
	$$ (-\Delta d_K)(x) \ge (-\Delta d_{\mathcal{P}})(x)  \ ,$$
since $ - \Delta f = \mathrm{Tr}(H(f)) $ where $H(f)$ stands for the Hessian of $f$. \\

One is then left to show that for any $x \in \HH^3$ and any geodesic hyperplane $\mathcal{P}$ such that $d(x, \mathcal{P}) \ge 1/2$ there is a constant $c > 0$ such that 
	$$ - \Delta d_{\mathcal{P}}(x) > 0 \ .$$
The proof consists in basic Riemannian geometry for which we outline the steps. Using the strong transitivity of the isometry group of $\HH^3$ and the Jacobi equation for variations of geodesics, one shows first that the hyperbolic metric of $\HH^3$ writes in exponential normal coordinates (with respect to $\mathcal{P}$) as
	$$ g_{\mathbb{H}^3} := dt^2 + \cosh(t)^2 d_{\mathbb{H}^2} \ ,$$
where $t = d_{\mathcal{P}}(x)$ and where $d_{\mathbb{H}^2}$ is the hyperbolic metric on the totally geodesic hyperplane $\mathcal{P}$. One then easily sees (using $\Delta f = \delta d f $) that 
$$ - \Delta t = 2 \tanh(t) \ge 2 \tanh(1/2) > 0 $$ 
for any $x \in \HH^3$ with $t = d_{\mathcal{P}}(x) \ge 1/2$, concluding. \hfill $\blacksquare$ $\blacksquare$  \\

Degenerate ends have been described more recently thanks to Minsky's work, see \cite{surveyminsky} and the references therein. The picture is now quite simple if $M_{\Gamma}$ is thick, they are roughly isometric to $\NN$.

\begin{theoreme}{ \cite[Corollary 14.8]{articlejuanbiringer}}
	\label{theoreme thick ends}
Let $\Gamma$ be a finitely generated Kleinian group such that $M_{\Gamma}$ is thick. Then $CC(M_{\Gamma})$ is roughly isometric to a graph which consists in gluing finitely many copies of $\NN$ on their common $0$.
\end{theoreme}

Note that in \cite{articlejuanbiringer} the authors do not mention a rough isometry but only a quasi-isometry: in this setting it is equivalent since we assumed $M_{\Gamma}$ to be thick. \\

Note also that in the case where $M_{\Gamma}$ is fully degenerate we have $CC(M_{\Gamma}) = M_{\Gamma}$. In particular, $M_{\Gamma}$ is roughly isometric to a graph consisting on finitely copies of $\NN$ glued together.

\subsection{Sullivan-Thurston's harmonic function.}

In this subsection we introduce the Sullivan-Thurston's harmonic function, which will ultimately plays the role of the weight we will use in order to estimate the large time behaviour of the heat kernel in the cases where the manifold $M_{\Gamma}$ carries both geometrically finite ends and fully degenerate ones.

\begin{theoreme}{\cite{articlesullivandemicylindre} \cite{articlebishopjones} \cite{articlejuanbiringer}}
	\label{theo harmonic function}
	Let $\Gamma$ be a finitely generated Kleinian group such that the manifold $M_{\Gamma}$ is thick and carries both degenerate and geometrically finite ends, then there exists a harmonic function $h_0 > 0$ such that 
		\begin{enumerate}
			\item the function $h_0$ has linear growth in degenerate ends: given a marked point $x_0 \in M_{\Gamma}$, there is two constants $a,b > 0$ satisfying that for all $x$ in the union of all the degenerate ends then
				$$ a^{-1} d(x_0,x) -b \le  h_0(x) \le a d(x_0,x) +b \ ;$$
			\item the harmonic function $h$ converges to $0$ in the geometrically finite ends.
		\end{enumerate}
\end{theoreme}

This theorem does not hold for fully degenerated manifold, they carry the Liouville property as a corollary of Ahfors' conjecture (now a theorem). \\

We will denote by $h_c$ the harmonic function $h_0 + c$, which enjoys the same properties than $h_0$, except that it converges to $c$ in the geometrically finite ends. \\

The story of this theorem is more intricate than the above statement might let think. In \cite[Theorem 3 page 135] {articlesullivandemicylindre}, carrying further Thurston's ideas, Sullivan showed that a non constant harmonic function has to have linear growth on manifolds roughly isometric to $\RR_+$ (that he called quasi-cylinder). Later on, Bishop and Jones \cite[lemma 1.7]{articlebishopjones} constructed a positive harmonic function in the full generality on any quotient of $\HH^3$ by a finitely generated Kleinian group, without requiring the ends to be quasi-cylinder. Moreover they classify them according to the number of degenerate ends. Roughly, they showed that the set of the harmonic functions described in Theorem \ref{theo harmonic function} is the positive convex cone generated by the finitely many harmonic functions $h_1,..., h_d$, where $d$ is the number of degenerate ends and $h_i$ is the only (up to a multiplicative constant) harmonic function which goes to infinity in the $i$-nth degenerate end and goes to $0$ in any other ends. The function $h_0$ appearing in Theorem \ref{theo harmonic function} can be taken as any combination of these functions with positive  coefficients. \\

Lacking a precise description of the degenerate ends, they could not manage to get the linear growth. Theorem \ref{theoreme thick ends} comes here to fill this lack by showing that the ends are actually roughly isometric to $\NN$ (or quasi-cylinders to follow Sullivan's terminology). 

\section{Proof of Theorem \ref{theo chaleur Kleinian}}

\label{sec proofs}
 
The two counting Theorems \ref{theo kleinian counting 1} and \ref{theo kleinian counting 2} announced in our introduction are corollaries of Theorems \ref{theo upper bound}, \ref{theo refined estimates} and \ref{theo chaleur Kleinian}. Namely, combined with Theorems \ref{theo upper bound} and \ref{theo refined estimates} one readily gets Theorem \ref{theo kleinian counting 1} and the first part of Theorem \ref{theo kleinian counting 2}. Second part of Theorem \ref{theo kleinian counting 2} follows from Theorem \ref{theo refined estimates} combined with second part of Theorem \ref{theo chaleur Kleinian} which gives that for every point $x \in M_{\Gamma}$ there is two constants $C_-$ and $C_+$ such that for $\rho$ large enough one has	
	$$ \frac{C_- \ e^{2\rho} }{\rho^{\frac{3}{2}}} \le N_{\Gamma}(x,x,\rho) \le \frac{C_+ \ e^{2\rho} }{\rho^{\frac{3}{2}}} \ . $$
This inequalities self improve immediately to the statement proposed in Theorem \ref{theo kleinian counting 2} since one has 
	$$ N_{\Gamma}(x,x, \rho - d(x,y)) \le N_{\Gamma}(x,y,\rho) \le N_{\Gamma}(x,x, \rho + d(x,y)) \ , $$ 
which gives the desired result, up to enlarging the constants. \\

It remains then to prove Theorem \ref{theo chaleur Kleinian}; estimating the large time behaviour of the heat kernel on thick degenerate hyperbolic manifolds.

\subsection{the fully degenerate case.}

\label{subsec fully degenerate} Recall that we want to prove the following 

\begin{theoreme}
	\label{theofullydeg}
	Let $M_{\Gamma}$ be fully degenerate, then there is two constants $C_-, C_+$ such that for any two points $x,y \in M_{\Gamma}$ there is $t_0 > 0$ such that for $t > t_0$ we have 
	$$ \frac{C_-}{\sqrt{t}} \le p_{\Gamma}(x,y,t) \le \frac{C_+}{\sqrt{t}} \ . $$
\end{theoreme} 

One can deduce the above theorem from \cite[Theorem 2.3 (i) page 10]{articlegrigosaloffcoste}. However, we prefer to give it here a proof in order to bypass all the refined technology involved in \cite{articlegrigosaloffcoste} which was developed for harder-to-get estimates. \\ 

In the case where the manifold $M_{\Gamma}$ carries only two degenerate ends, as it is for the historical examples given by degeneration of quasi-Fuchsian manifolds \cite{livremcmumu}, it follows from Theorem \ref{theoreme thick ends} that $M_{\Gamma}$ is roughly isometric to $\mathbb{Z}$. If one would be able to compare the Brownian motion on $M_{\Gamma}$ and the random walk on $\mathbb{Z}$, the intuition behind Theorem \ref{theofullydeg} becomes clear since the local limit theorem for a random walk on $\mathbb{Z}$ is well known to behave as in the above theorem.  \\

\begin{figure}[h!]
\begin{center}
	\def\svgwidth{0.8 \columnwidth}
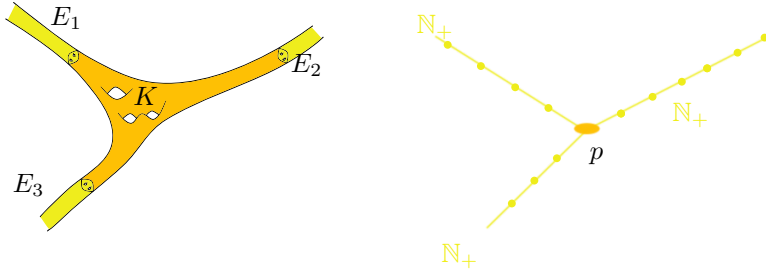
\caption{On the left, a fully degenerate manifold with 3 ends which is roughly isometric to the graph $G_3$ on the right.}
\label{figurefulldeggraph}
	\end{center}
\end{figure}

\textbf{Proof.} We want to use Theorem \ref{theo grigo equivalence}. In order to get the conclusion of Theorem \ref{theo grigo equivalence}, one must verify that $M_{\Gamma}$ is roughly isometric to a graph which carries the doubling property and the Poincaré inequality. Recall that the convex core of a fully degenerate group is equal to the whole manifold. Therefore, by Theorem \ref{theoreme thick ends}, a thick manifold $M_{\Gamma}$ associated to a fully degenerate Kleinian group is roughly isometric to the graph built from a rooted point on which we attach $d$ copies of $\mathbb{N}$ - corresponding to the degenerate ends - as in Figure \ref{figurefulldeggraph}. We denote such a graph by $G_d$. The doubling volume property is immediate for the graph $G_d$ since we have for all $r > 0$ 
	$$  r \le \vol\big(B(x,r)\big) \le d \ r \ .  $$

The proof proposed below of the Poincaré inequality for the graphs $G_d$ follows the same line than the short - less than a page - one proposed by Kleiner and Saloff-Costes in \cite[Theorem 2.2]{artkleiner} at the very beginning of the paper. Note however that the setting is different since  \cite{artkleiner} deals with Cayley graphs of groups. \\
	 	 
Here $\mu$ stands for the counting measure on the graph $G_d$. We shall use the integral notation even if the support of the measure $\mu$ is discrete: we reserve the use of the summation notation along the proof to emphasis the key step. \\

Let $ x \in G_d$, $r > 0$ and $f$ be a function on $G_d$, we start off the left member of Equation \eqref{eqpoincare}:	 
	 $$ \inte{B(x,r)} ( f- f_r )^2 d \mu =
	  \inte{B(x,r)} \left( \frac{1}{ \vol\big(B(x,r)\big)}   \inte{B(x,r)} (f(z) - f(y)) \ d \mu(y) \right)^2  d \mu(z) \ .$$
Since we have the lower bound $\vol\big(B(x,r)\big) \ge r$ we get
$$ \inte{B(x,r)} ( f- f_r )^2 d \mu \le \frac{1}{r^2}
	  \inte{B(x,r)} \left( \inte{B(x,r)} f(z) - f(y) \ d \mu(y) \right)^2  d \mu(z) \ . $$
	  
Using Cauchy-Schwartz' inequality one has
 $$ \left( \inte{B(x,r)} f(z) - f(y) \ d \mu(y) \right)^2  d \mu(z) \le \vol \big( B(x,r) \big) \inte{B(x,r)}  \left( f(z) - f(y) \right)^2  \ d \mu(y) \ . $$
 
We now use the upper bound $\vol \big( B(x,r) \big) \le d \cdot r$ to get 
\begin{equation}
	\label{eqprooffulldege}
	 \inte{B(x,r)} ( f- f_r )^2 d \mu \le \frac{d}{r}  \inte{B(x,r)} \inte{B(x,r)}  \left( f(z) - f(y) \right)^2  \ d \mu(y) d \mu(z) \ .
\end{equation}

Let us denote by $\gamma_{y,z} = (\omega_0 = y,\omega_1,...,\omega_n = z)$ the (unique) geodesic of $G_d$ relating the points $x$ and $y$, where $n = d(x,y)$. Writing
	$$ f(z) - f(y) \le \somme{1 \le k \le n} \big (f(\omega_{i}) - f(\omega_{i-1}) \big) \ ,$$
and using again the Cauchy Schwartz identity one gets  
  $$ ( f(z) - f(y) )^2 \le n \somme{1 \le i \le n} \big (f(\omega_{i}) - f(\omega_{i-1}) \big)^2 \ . $$\\
 Note that $n \le 2r$ since $x,y \in B(x,r)$. Note also that
$$ \somme{1 \le k \le n} \big (f(\omega_{i}) - f(\omega_{i-1}) \big)^2 \le \somme{\omega \in \gamma_{y,z} } (\delta f )^2 ( \omega) \ . $$

Recall that the discrete derivative operator $\delta$ is defined in Section \ref{sec weighted heat kernel}. Which gives 
 $$ ( f(z) - f(y) )^2 \le 2 r \somme{\omega \in \gamma_{y,z} } (\delta f )^2 ( \omega) \ .$$
 
Looking backward at \ref{eqprooffulldege}, we have 
	$$	 \inte{B(x,r)} ( f- f_r )^2 d \mu \le 2 d  \inte{B(x,r)} \inte{B(x,r)}  \somme{\omega \in \gamma_{y,z} } (\delta f )^2 ( \omega)   \ d \mu(y) d \mu(z) \ . $$

Any geodesic $\gamma_{z,y}$ of $G_d$ is exactly given by its endpoints, and, again because $\gamma_{z,y}$ is a geodesic, a point $p \in B(x,r)$ belongs to $\gamma_{z,y}$ at most once. Since we have prescribed in the above summation the endpoints of the geodesics to lie in $B(x,r)$, there is at most $(d \ r)^2$ different such geodesics and thus at most $(d \ r)^2$ occurrences of the term $( \delta f)^2(p)$. Therefore,
\begin{equation}	
	\label{eqrpoincarelineargrowth}
	 \inte{B(x,r)} ( f- f_r )^2 d \mu  \le 2 d^3 \ r^2 \inte{ B(x,r)} (\delta f)^2 (p) \ d \mu (p) \ ,
\end{equation}
	 which concludes the proof setting $P = 2 \ d^3$. \hfill $\blacksquare$ \\

Note that the proof could be adapted to the case of a graph of linear growth which would give 
\begin{align*}	
	 \inte{B(x,r)} ( f- f_r )^2 d \mu & \le P \ r^2 \inte{ B(x,2r)} (\delta f)^2 (p) \ d \mu (p) \ ,
\end{align*}	 	
for some $C$, which is exactly what is required for carrying a Poincaré inequality. Indeed, with Inequality \eqref{eqrpoincarelineargrowth} we obtained something stronger than what required since the right integral is over $B(x,r)$. 
\subsection{The mixed type case.} Recall here the statement dealing with the mixed type case.

\begin{theoreme}
\label{theo chaleur mixed type}
Let $M_\Gamma$ be a mixed type hyperbolic manifold then,
\begin{itemize}
	 \item  for every $x,y \in M_{\Gamma}$ there is a constant $C_+$ such that for $t$ large enough
			$$ p_{\Gamma}(x,y,t) \le \frac{C_+}{t^{\frac{3}{2}}} \ . $$
	\item Moreover, for any $x \in M_{\Gamma}$ there is a constant $C_-$ such that
		$$ \frac{C_-}{t^{\frac{3}{2}}} \le p_{\Gamma}(x,x,t) $$
\end{itemize}	 
\end{theoreme}

The end of this paper is devoted to the proof of this theorem. The intuition behind Theorem \ref{theo chaleur mixed type} is that whenever the Brownian motion enters the geometrically finite end, which looks roughly like a hyperbolic space, it has to choose a point at infinity and goes toward it. Therefore, it tends to escape every compact set; it seems then natural to compare, in large time, the behaviour of the heat kernel to the one of a Brownian motion defined on the degenerate ends stopped whenever it goes out the convex core. A discrete analogue would be a random walk defined on $\mathbb{N}$ - the degenerate end and stopped at $0$, known to have a return probability of the order of $n^{-\frac{3}{2}}$. \\

To get the two-sided estimate, we use Theorem \ref{theo lower bound} which needs an upper bound, so that we start in proving the upper bound.

\subsection{The upper bound.}

\begin{figure}[h!]
\begin{center}
	\def\svgwidth{0.6 \columnwidth}
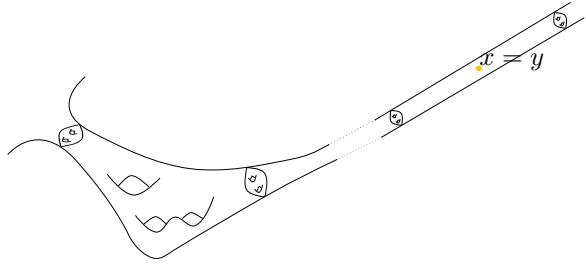
\caption{One expects the heat kernel $p_{\Gamma}(x,y,t)$, taken at $x=y$ far away in the degenerate end, to behave as the heat kernel of $\RR$ to which it is roughly isometric. This kernel behaves like $\frac{1}{\sqrt{t}}$, which prevents one to show either that $M_{\Gamma}$ carries a $S_{2,6}$ or the parabolic Harnack inequality. }
\label{figuremixedtype}
	\end{center}
\end{figure}

In order to show that the upper bound of Theorem \ref{theo chaleur mixed type} holds, one would be tempted to show that $M_{\Gamma}$ directly carries a Sobolev inequality $S_{2,6}$, meaning to show that there is a constant $C > 0$ such that for any smooth compactly supported function $f$ of $M$ we have
	$$ || f ||_{L^6(\mu_g)} \le C || \nabla f ||_{L^2(\mu_g)} $$
 which is known to imply the desired upper bound, see for example \cite[Corollary 14.23 page 383]{bookgrigorheatkernel}. Unfortunately, such a direct approach will surely fail since it gives an uniform upper bound on the heat kernel, regardless of where are the points $x$ and $y$, see Figure \ref{figuremixedtype}. Note that, for the same reason, there is no chance to prove that $M_{\Gamma}$ carries the parabolic Harnack inequality. \\
 
 In order to circumvent this difficulty, we weight $M_{\Gamma}$ according to the square of the Sullivan-Thurston harmonic function given by Theorem \ref{theo harmonic function}. \\

Recall that we denoted by $h_1$ the Thurston's Sullivan harmonic function which tends to $1$ in the geometrically finite end and that we denoted by $p_{h_1^2}(x,y,t)$ the heat kernel associated to the weighted Riemannian manifold $(M_{\Gamma},g, \mu_{h_1^2})$. Let us  come to the main proposition.

\begin{proposition}
	\label{propupperboundheatkernel}
		There are constant $C, t_0 > 0$ such that for any $t > t_0$ we have for any $x,y \in M_{\Gamma}$ 
	$$ p_{h_1^2}(x,y,t) \le C \ t^{-\frac{3}{2}} \ . $$
\end{proposition}

Note that this readily gives the desired upper bound for the usual heat kernel since, because of Proposition \ref{proposition heat kernels related}, we have 
 $$ p_{h_1^2}(x,y,t) = \frac{ p_{\Gamma}(x,y,t)}{h_1(x) \ h_1(y)} \ .$$ 

We rely on Theorem \ref{theocomparison} to prove the above Proposition. One should start to verify the
\begin{lemma}
	\label{lemmalocal}
	The weighted Riemannian manifold $(M_{\Gamma},g, \mu_{h_1^2})$ satisfies the $(P_2)_0$ and the $(DV)_0$ properties.
\end{lemma}

This done, we will conclude by using Theorem \ref{theocomparison} together with the

\begin{proposition}
	\label{propmodel}
	The weighted Riemannian manifold $(M_{\Gamma},g, \mu_{h_1^2})$ is roughly isometric to a graph $G$ which satisfies a Sobolev inequality $S_{2,6}$.
\end{proposition}

It remains then to prove the two above statements. \\

\textbf{Proof of lemma \ref{lemmalocal}.} We introduce some notations. If $\sigma$ is a positive function on a Riemannian manifold $(M,g)$, we denote by

\begin{itemize}
	\item $\mathcal{V}_g(x,r)$ the $\mu_{g}$-volume of the ball $B(x,r)$;
	\item $\mathcal{V}_\sigma(x,r)$ the $\mu_{\sigma}$-volume of the ball $B(x,r)$;
	\item $\sigma_-(x,r)$ the infimum of $\sigma$ over the ball $B(x,r)$; 
	\item $\sigma_+(x,r)$ the suppremum of $\sigma$ over the ball $B(x,r)$.
\end{itemize}

Note that for any $x$ the function $ r \mapsto \sigma_+(x,r)$ is increasing and $ r \mapsto\sigma_-(x,r)$ is decreasing. \\

We start with the easiest part, the property $(DV)_0$. We fix $r_0 > 0$ and $r \le r_0$. With the above notation we have 
$$  \mathcal{V}_{\sigma}(x,2r) \le \sigma_+(x,2r) \ \mathcal{V}_g(x,2r), $$
But, since $(M_{\Gamma},g)$ is hyperbolic with positive injectivity radius, it satisfies the property $(DV)_0$ so that there is a constant $C_{r_0}$ such that 
$$  \mathcal{V}_g(x,2r) \le C_{r_0} \mathcal{V}_g(x,r). $$

In particular,
\begin{align*}
	\mathcal{V}_{\sigma}(x,2r) & \le C_{r_0} \sigma_+(x,2r) \ \mathcal{V}_{g}(x,r) \\
	& \le C_{r_0} \frac{\sigma_+(x,2r)}{\sigma_-(x,r)}  \mathcal{V}_{\sigma}(x,r) \\
	& \le C_{r_0} \frac{\sigma_+(x,2r_0)}{\sigma_-(x,2r_0)}  \mathcal{V}_{\sigma}(x,r) \ ,
\end{align*}
since $r \le r_0$. In our setting we have $\sigma = h_1^2$, therefore:
$$ \frac{\sigma_+(x,2r_0)}{\sigma_-(x,2r_0)} = \frac{({h_1^2})_+(x,2 r_0)}{({h_1^2})_-(x,2 r_0)}  =  \left( \frac{{h_1}_+(x,2 r_0)}{{h_1}_-(x,2 r_0)} \right)^2 \ . $$ 

Using Harnack's inequality ($h_1$ is harmonic), the right member of the above inequality is uniformly bounded by a constant which only depends on $r_0$, which concludes for what concerns the $(DV)_0$ property. \\

Let us come to the condition $(P_2)_0$. We start with noticing that for any constant $c \in \RR$ we have 
	$$ V_{\sigma}(f,x,r) := \inte{B(x,r)} \left( f -\frac{1}{\mathcal{V}_{\sigma}(x,r)} \inte{B(x,r)} f d \mu_{\sigma} \right)^2 d \mu_{\sigma} \le \inte{B(x,r)} \left( f -c \right)^2 d \mu_{\sigma} \ .  $$

In particular, setting 
	$$ c := \frac{1}{\mathcal{V}_{g}(x,r)} \inte{B(x,r)} f d \mu_{g} $$ 
we get 
	\begin{align*}
			V_{\sigma}(f,x,r) & \le \inte{B(x,r)} \left( f - \frac{1}{\mathcal{V}_{g}(x,r)} \inte{B(x,r)} f d \mu_{g}  \right)^2 d \mu_{\sigma} \\
			 & \le \sigma_+(x,r) \inte{B(x,r)} \left( f - \frac{1}{\mathcal{V}_{g}(x,r)} \inte{B(x,r)} f d \mu_{g}  \right)^2 d \mu_{g} \ .
	\end{align*}

The local Poincaré inequality for the non weighted hyperbolic manifold $(M_{\Gamma},g)$ gives a constant $C_{r_0}$ such that for any smooth function $f : M_{\Gamma} \to \RR$ we have  
	$$ \inte{B(x,r)} \left( f -\frac{1}{\mathcal{V}_{g}(x,r)} \inte{B(x,r)} f d \mu_{g} \right)^2 d \mu_{g} \le \ C_{r_0} \ r^2 \inte{B(x,2r)} | \nabla f|^2 d\mu_{g} \ . $$
And then 
	\begin{align*}
		V_{\sigma}(f,x,r) & \le \sigma_+(x,r) \ C_{r_0} \ r^2 \ \inte{B(x,2r)} | \nabla f|^2 d\mu_{g} \\ 
				 & \le  \frac{\sigma_+(x,r)}{\sigma_-(x,2r)} \ C_{r_0} \ r^2 \  \inte{B(x,2r)} | \nabla f|^2 d\mu_{\sigma} \\
				 & \le \frac{\sigma_+(x,2r_0)}{\sigma_-(x,2r_0)} \ C_{r_0} \ r^2 \  \inte{B(x,2r)} | \nabla f|^2 d\mu_{\sigma} \ ,
	\end{align*}
since $r \le r_0$. We conclude as we did for property $(DV)_0$ by using Harnack's inequality. \hfill $\blacksquare$ \\

Let us come now come to the global part, the proof of Proposition \ref{propmodel}. \\

\textbf{Proof of Proposition \ref{propmodel}.} The manifold $M_{\Gamma}$ being given as a codimension $0$ compact submanifold $K$ with boundary on which we attached finitely many ends, it is roughly isometric to discretisations of these finitely many ends glued together at some definite point. We denote by $GF_1,... GF_q$ the geometrically finite ends and by $D_1,...D_d$ the degenerate ones. \\

From Theorem \ref{theoreme thick ends}, the convex core of $M_{\Gamma}$ (which corresponds to the degenerate ends together with $K$) is roughly isometric to the graph $G_d$ (appearing in Proposition \ref{theofullydeg}). Recall that $G_d$ was defined as the graph obtained in gluing together $d$ different copies of $\NN$ on their common zero denoted $x_0$. It remains then to glue the discretisation of the geometrically finite ends to $x_0$ (we will not use anything about them except that they have positive bottom of the spectrum). We then recover a graph $G_{d,p}$ which is roughly isometric to the hyperbolic manifold $(M_{\Gamma},g)$. \\

In particular, because of Theorem \ref{theo harmonic function}, the weighted Riemannian manifold $(M_{\Gamma},g, h_1^2)$ is roughly isometric to the weighted graph $(G_{d,p}, m^2)$ where $m : G_{d,p} \mapsto  \NN$ is the constant function equal to one in the geometrically finite ends and the identity on any degenerate ends (the weight of a vertex $x \in \NN \setminus \{ 0\}$ is exactly given by the number $x$). Note that the discretisations of geometrically finite ends are not weighted. \\

Let us now prove that the weighted graph $(G_{d,p}, m^2)$ carries a $S_{2,6}$. We shall actually see that it is the case for every end independently and we will conclude using the lemma below. The analogous for Riemannian manifold of the following lemma was proved in \cite{articlecarron}, but the proof goes the same way (it is even simpler) for weighted graphs and is left to the reader. 

\begin{lemma}{ \cite[Proposition 2.5]{articlecarron}}
	Let $N \in \NN_+$, $p > 2$ and $(G_1,x_1,m_1)$... \\ $(G_m,x_N,m_N)$ be finitely many weighted rooted graphs carrying a $S_{2,p}$. 
Then the weighted graph obtained in gluing all the $G_i$s along their roots carries a $S_{2,p}$.
\end{lemma}

We are left to show that the discretisations of both geometrically finite ends and weighted degenerate ones satisfies a $S_{2,6}$. \\

 We start with the discretisation of a geometrically finite end. Note first that the weight $m$ on $G_{d,p}$ is constant equal to one on the discretisations of the geometrically finite ends. \\ 

Recall that geometrically finite ends have positive bottom of the spectrum by Proposition \ref{propspectralgap}. Note that, by definition, having positive bottom of the spectrum is equivalent to carrying a Sobolev inequality $S_{2,2}$. Since carrying an $S_{2,2}$ on a manifold is equivalent to carrying one on any of its discretisation, the discretised of any geometrically finite end also carries a $S_{2,2}$ (see \cite[Proposition 6.5]{artcoulhonsaloffisometrie}). \\

One is then left to show that if a metric graph $X$ satisfies a $S_{2,2}$ then it also satisfies a $S_{2,6}$. This follows readily from the fact that the embedding from $L^2(X,\mu)$ to $L^p(X,\mu)$ is continuous for any $p \ge 2$ (in particular for $p =6$). Indeed, such a continuous embedding holds whenever the measure $\mu$ gives uniform positive mass to any non zero measure set, as it the case for the counting measure on a graph. 

\begin{remark} Note that the reference actually asserts that if $(\RR_{\ge 1}, dx^2, x^2 dx)$ carries a $S^{\infty}_{2,6}$ then $(\NN_{\ge 1},n^2)$ carries a $S_{2,6}$. But carrying a $S_{2,6}$ for $(\RR_{\ge 1}, dx^2, x^2 dx)$ is stronger than carrying a $S^{\infty}_{2,6}$, since the last is a version 'localised at infinity' of the previous one, see \cite[end of the page 705]{artcoulhonsaloffisometrie}. 
\end{remark}

We will use the fact that $\RR^3$ carries a $S_{2,6}$: there is a constant $C > 0$ such that for any compactly supported smooth function $F : \RR^3 \to \RR$ 
	$$ || F||_{L^6(\RR^3)} \le C  || \nabla F||_{L^2(\RR^3)} \ . $$
On the one hand, if $F$ is taken radial ($F(x) = f(r)$ where $r = ||x||$) we have
	$$ || F||_{L^6(\RR^3)} = (4 \pi)^{1/6} || f ||_{L^6(\RR_+, r^2 dr)} \ ,$$ 
since 
	$$ \inte{\RR^3} \ F^6 \ d \mu_{\RR^3} = \inte{\RR_+} \ f^6(r) \ 4 \pi \ r^2 dr \ . $$
On the other hand, we have $$ || \nabla F||_{L^2(\RR^3)} = 2 \sqrt{\pi} \ || f' ||_{L^2(\RR_+, r^2 dr)} \ . $$
In particular there is a constant $C > 0$ such that for any compactly supported function on $\RR_{\ge 1}$ we have
 $$  (4 \pi)^{1/6} \ || f ||_{L^6(\RR_{\ge 1}, r^2 dr)} \le 2 \ C \ \sqrt{\pi} \ || f' ||_{L^2(\RR_{\ge 1}, r^2 dr)} \ .$$

In other words, $(\RR_{\ge 1}, dx^2, x^2 dx)$ carries a $S_{2,6}$.  \hfill $\blacksquare$

\subsection{The lower bound} We conclude this article by showing the lower bound in the mixed type case: for any $x \in M_{\Gamma}$ there is $t_0 \ge 0$ and a constant $C > 0$ such that for any $t > t_0$ we have 
	$$ p_{\Gamma}(x,x,t) \ge  \frac{C}{t^{\frac{3}{2}}} \ . $$

In fact, we will prove that the above lower bound holds for $t_0 = 0$. Note first, using Proposition \ref{proposition heat kernels related} again, that one only needs to show that for any $x \in M_{\Gamma}$ there is a constant $C > 0$ such that for any $t > 0$ we have 
\begin{equation}
	\label{eqlowerbound}
	 p_{h_0^2}(x,x,t) \ge  \frac{C}{t^{\frac{3}{2}}} \ .
\end{equation}
Recall that $h_0$ is the Thurston Sullivan harmonic function which is asymptotically $0$ in the geometrically finite ends. \\

In order to get \eqref{eqlowerbound}, we now want to use Theorem \ref{theo lower bound}. To use this theorem one much check both that for a given $x \in M_{\Gamma}$ and for any $t > 0$ there is a constant $C_+$ such that 
	$$ p_{h_0^2}(x,x,t) \le C_+ \ \mu_{h_0^2} \Big(B(x,\sqrt{t}) \Big)^{-1} \ , $$
and that $(M_{\Gamma},g, d \mu_{h_0^2})$ is doubling at $x$. \\

Because of Proposition \ref{propupperboundheatkernel} and since, using again \ref{proposition heat kernels related}, we already now that for any $x \in M_{\Gamma}$ there are constants $C, t_0 > 0$ such that for $t > t_0$ we have 	
	$$	p_{h_0^2}(x,x,t) \le  \frac{C}{t^{\frac{3}{2}}} \ .$$ 

In fact, the same upper bound also holds and for any $0 < t \le t_0 $, $x$ being fixed. This is immediate since we have by Proposition \ref{proposition heat kernels related} and \ref{lemma heat kernel quotient} (with the notation of \ref{lemma heat kernel quotient}) 
	$$ p_{h_0^2}(x,x,t) = \frac{1}{h_0(x)^2} \ \somme{\gamma \in \Gamma} \ p_3(\tilde{x}, \gamma \cdot \tilde{x}, t) \ . $$
	
Indeed, in small time, the only term contributing to this summation corresponds to the identity of $\Gamma$, giving the desired upper bound using the explicit formula given by \ref{theo heat kernel 3-hyperbolic}. To sum up, we already had proven that for any $x \in M_{\Gamma}$ and for $t > 0$ 
		$$	p_{h_0^2}(x,x,t) \le  \frac{C}{t^{\frac{3}{2}}} \ .$$ 

The following geometric proposition shows in particular that 
	$$\mu_{h_0^2} \big(B(x,\sqrt{t}) \big) \le C t^{\frac{3}{2}} $$ 
and that $(M_{\Gamma},g, d \mu_{h_0^2})$ carries the doubling property around $x$. One can then apply Theorem \ref{theo lower bound} to get the desired lower bound and conclude. Let us then show

\begin{proposition}
	\label{prop doubling}
	For any $x \in M_{\Gamma}$ there is two constant $C_+,C_-$ such that for any $t \ge 0$
	 $$  C_- t^{3} \le \mu_{h_0^2} \Big(B(x,t) \Big) \le C_+ t^{3} \ . $$ 
	 
In particular the manifold $M_{\Gamma}$ enjoys the volume doubling property around $x$.
\end{proposition}

\textbf{Proof.} Since we know the harmonic function to have linear growth in any degenerate ends we readily gets there are constants $C_-, C_+$ such that 
	 $$  C_- t^{3} \le \mu_{h_0^2} \Big( CC(M_{\Gamma}) \cap B(x,t) \Big) \le C_+ t^{3} \ . $$ 
 which gives that the volume of the convex core is of order $t^{3}$. Recall that $CC(M_{\Gamma})$ is the convex core of $\Gamma$ and contains all the degenerate ends. \\

One would be then able to conclude if the $\mu_{h_0^2}$ volume of the complementary of $CC(M_{\Gamma})$, the geometrically finite ends, is finite. That is to say that the harmonic function $h_0$ is in $L^2(M \setminus CC(M_{\Gamma}), \mu_g)$. This readily follows from the

\begin{lemma}
	Let $(M,g)$ be a hyperbolic manifold and $E$ an end of $M$. If $E$ has infinite volume, positive injectivity radius and positive bottom of the spectrum, then any positive harmonic function $h$ on $E$ such that $h(y) \tends{y \to \infty} 0$ belongs to $L^2(M,g)$.
\end{lemma}

\textbf{Proof.} To prove the above lemma we shall use the notion of barrier functions. Let $E$ be an end of a Riemannian manifold $(M,g)$ manifold such that $\mu_g(E)$ is infinite and $\lambda_0(E) > 0$. Let $p$ be a fixed point in $M$, for any $R > 0$ denote by $E(R) := B(p,R) \cap E$. So that for $R$ large enough one has $\partial E(R) = \partial E \cup (\partial B(p,R) \cap E) $. Let $h_R$ be the solution of the following Dirichlet problem
\begin{equation*}	
	\left\{
		\begin{array}{l}
			\Delta h_R = 0 \\
			h_R = 1 \ \text{ on } \ \partial E  \\
			h_R = 0 \ \text{ on } \ \partial B(p,R) \cap E
		\end{array} 
	\right. 
\end{equation*}
Given a sequence of $R_n  \tends{n \to \infty} \infty$, there is a subsequence of $(h_{R_n})_{n \in \NN}$ which converges uniformly on any compact to a positive harmonic function $h_{\infty}$ such that $h_{\infty} = 1$ on $\partial E$. Any such a limit is called a barrier function, see \cite[section 20]{bookligeometricanalysis} for existence. \\

We will prove first that there is a barrier function which is in $L^2(E,g)$ and we will conclude using the maximum principle to show that the harmonic function $h_0$ is also in $L^2(E,g)$. \\

Since $E$ has infinite volume and $\lambda_0(E) > 0$ the end $E$ must carry a non constant barrier function  \cite[Proposition 22.5]{bookligeometricanalysis} (see \cite[Section 20]{bookligeometricanalysis} for the definition of a non parabolic end). Up to subtracting a constant function and multiplying by a positive constant one can assume that $\limiinf{y \to \infty} \ h_{\infty} = 0$. Now, applying \cite[Corollary 22.3]{bookligeometricanalysis} (see the remark in the end of the proof) to such a function gives that $h_{\infty} \in L^2(E,g)$. In particular, $h_{\infty}(y) \tends{y \to \infty} 0$ by harnack's inequality. \\

As already emphasised, we conclude by the maximum principle. Let $R_n \to \infty$ such that $ h_{R_n} \to h_{\infty}$ with $h_{\infty}(y) \tends{y \to \infty} 0$ and $h_{\infty} \in L^2(E,g)$. By the maximum principle one has for any $n \in \NN$:
	$$ h_0 \le \left( \supr{\partial E} h_0\right) \ h_{R_n} + \supr{\partial B(p,R_n) \cap E} \ h_0 \ .$$
Since $h_0 \tends{y \to \infty} 0$ in restriction on $E$ by Theorem \ref{theo harmonic function} we have, letting $n \to \infty$,
	$$ h_0 \le \left(\supr{\partial E} \  h_0 \right) \ h_{R_{\infty}} \ , $$
concluding since $h_0 > 0$.  \hfill $\blacksquare$ $\blacksquare$

\bibliographystyle{plain}
\bibliography{bibliography} 

\end{document}